\documentclass[a4paper,10pt,twocolumn]{article}
\pdfoutput=1

\usepackage{graphicx}
\usepackage{grffile}
\usepackage{xcolor}
\usepackage{titling} 
\usepackage{soul}    

\usepackage{url}
\usepackage[pdftex,colorlinks=true,linkcolor=blue,citecolor=blue,urlcolor=blue]{hyperref}

\usepackage[expproduct=cdot]{siunitx} 
\usepackage{enumitem} 
\usepackage{epstopdf}

\usepackage{relsize} 
\usepackage{amsmath}
\usepackage{amssymb}\renewcommand{\geq}{\geqslant}\renewcommand{\leq}{\leqslant}
\usepackage{bbm} 
\usepackage[utopia,greekuppercase=italicized]{mathdesign}\edef\partial{\mathchar\number\partial\noexpand\!} 

\usepackage[T1]{fontenc}

\usepackage{ctable} 
\usepackage{tabularx}

\usepackage[numbers,compress]{natbib}

\usepackage{subfigure}

\definecolor{royalblue4}{HTML}{27408B} 
\definecolor{red4}{HTML}{8B0000} 
\definecolor{green4}{HTML}{008b00}

\usepackage{dblfloatfix}  
\usepackage{fixltx2e}  

\usepackage[font=footnotesize,labelfont=bf,labelsep=endash,margin=0pt]{caption} 
\usepackage[title,header]{appendix}

\newlength{\myleftmargin} 
\newlength{\mytopmargin} 
\newlength{\myrightmargin} 
\newlength{\mybottommargin} 
\usepackage[runin]{abstract}

\let\paragraphold\paragraph
\renewcommand*{\paragraph}[1]{\paragraphold{#1.}} 
\newcommand{\keywords}[1]{\vspace{2mm}\noindent\textbf{Keywords:} #1} 
\newcommand{\pagewidetitle}[3] 
{
    \twocolumn
        [
            \vskip-5mm
            \begin{@twocolumnfalse}%
                #1
                #2
                \vspace{5mm}
            \end{@twocolumnfalse}%
        ]
        #3
}

\newlength{\figurewidth}

\definecolor{pbuenzli}{HTML}{00b7eb} 
\definecolor{pbuenzlinew}{HTML}{008b00}
\definecolor{pbuen}{HTML}{A020F0}

\makeatletter
\long\def\tlist@if@empty@nTF #1{
\expandafter\ifx\expandafter\\\detokenize{#1}\\
\expandafter\@firstoftwo
\else
\expandafter\@secondoftwo
\fi
}

\newcommand*\pbuenzli[2][]{
{\color{pbuenzli}#2}
\tlist@if@empty@nTF{#1}{}{\footnote{\color{pbuenzli}\it #1}}
}
\makeatother

\newcommand{\ups}{{S}}

\renewcommand{\d}{\mathrm{d}}
\newcommand{\p}{\partial}

\newcommand{\e}{\mathrm{e}}

\renewcommand{\b}[1]{{\boldsymbol{#1}}} 

\newcommand{\Tr}{\mathrm{Tr}}

\usepackage{relsize}

\newcommand{\kform}{\text{$k$}}

\begin{document}



    \title{\bf A level-set method for the evolution of cells and tissue during curvature-controlled growth}




	\renewcommand{\thefootnote}{\fnsymbol{footnote}}%
	\author{Mohd Almie Alias$^\text{a}$\footnotemark[1], Pascal R Buenzli$^\text{b}$}
	
	\date{\normalsize \vspace{-2mm}
	$^\text{a}$ Center for Modelling and Data Science, Faculty of Science and Technology, Universiti Kebangsaan Malaysia,\\ 43600 Bangi, Selangor D. Ehsan, Malaysia\\$^\text{b}$School of Mathematical Sciences, Queensland University of Technology, Brisbane QLD 4001, Australia\\\vskip 1mm \normalsize	
	\today\vspace*{-5mm}}
	
	\pagewidetitle{
		\maketitle
	}{
      \begin{abstract}
Most biological tissues grow by the synthesis of new material close to the tissue's interface, where spatial interactions can exert strong geometric influences on the local rate of growth. These geometric influences may be mechanistic, or cell behavioural in nature. The control of geometry on tissue growth has been evidenced in many in-vivo and in-vitro experiments, including bone remodelling, wound healing, and tissue engineering scaffolds. In this paper, we propose a generalisation of a mathematical model that captures the mechanistic influence of curvature on the joint evolution of cell density and tissue shape during tissue growth. This generalisation allows us to simulate abrupt topological changes such as tissue fragmentation and tissue fusion, as well as three dimensional cases, through a level-set-based method. The level-set method developed introduces another Eulerian field than the level-set function. This additional field represents the surface density of tissue synthesising cells, anticipated at future locations of the interface. Numerical tests performed with this level-set-based method show that numerical conservation of cells is a good indicator of simulation accuracy, particularly when cusps develop in the tissue's interface. We apply this new model to several situations of curvature-controlled tissue evolutions that include fragmentation and fusion.

	\end{abstract}
}{
\protect\footnotetext[1]{Corresponding author. \texttt{mohdalmie@ukm.edu.my}}%
\renewcommand{\thefootnote}{\arabic{footnote}}%
}

\section{Introduction} \label{Introduction}
Biological tissues grow, change shape and material properties under strong geometric controls~\cite{Nelson2005,Rumpler2008,Bidan2012,Bidan2013,Jin2018}. Some of these controls are mechanistically induced by the tissue's evolving geometry, such as the crowding and spreading of tissue constituents due to spatial constraints~\cite{Alias2017}, while other geometric controls relate to influences on cell behaviours such as level of activity, proliferation, and death~\cite{Alias2018}. Mathematical models of tissue growth and morphogenesis commonly model the overall geometric control by mean curvature flows, whereby the normal velocity of the tissue interface is simply proportional to the local curvature~\cite{Rumpler2008,Bidan2012,Guyot2014,Sanaei2019}. These phenomenological models do not consider the cellular component of tissue synthesis, and as such, are unable to disentangle the mechanistic and cell behavioural origins of curvature dependence in tissue growth. Yet, this distinction is important for the correct interpretation of experimental data. The crowding and spreading of tissue constituents induced by curvature is an unavoidable effect that must be factored out to analyse cell behaviours~\cite{Alias2018}.

In previous works, we have proposed a cell-based mathematical model of tissue growth to account for the mechanistic influence of curvature on cell crowding and cell spreading in the co-evolution of cell density and tissue interface~\cite{Alias2017,Alias2018}. This mathematical model reduces to a specific type of hyperbolic curvature flow, in which the normal acceleration of the interface is proportional to curvature~\cite{LeFloch2008,Dexing2009,He2009,Ishida2017,Grayson1987}. The hyperbolic character of this curvature flow gives rise to a rich set of interface movement patterns. Depending on the amount of lateral diffusive damping (related to cell migratory behaviours), these movement patterns include oscillatory interface motion, emergence of cusps propagating sideways as shock waves, efficient smoothing of initial interface irregularities, and shock and rarefaction waves emerging at concavities and convexities~\cite{Alias2017,Alias2018}. However, these previous works are restricted in two important ways: (i) they do not consider abrupt topological changes of the interface that may occur when separate regions of the tissue merge, or when the tissue becomes fragmented into disconnected pieces; (ii) the model is two-dimensional and its extension to three dimensional space is nontrivial. In biology, topological changes commonly arise in morphogenesis and tissue growth. They occur for example in tumour growth (fusion of tumour nodules or of tumour fingers)~\cite{Macklin2006,Lowengrub2010,Wise2008}, wound healing~\cite{Poujade2007,Vermolen2007}, tissue involution processes, tissue engineering bioscaffolds~\cite{Guyot2014,Guyot2015,Guyot2016,Paris2017}, and in bone consolidation and fragmentation~\cite{Maggiano2016b,Bell2001,Kinney1998,Jordan2000}.

The level set method is a common and successful technique to numerically evolve interfaces undergoing topological changes~\cite{Sethian1999,Osher2003,Gibou2018}. It has been used to describe the evolution of biological tissues in several instances before~\cite{Macklin2003,Macklin2006,Wise2008,Guyot2014,Guyot2015,Guyot2016,Vermolen2007}. However, in these studies, the population of tissue-synthesising cells is not considered, and the tissue's interface velocity is usually assumed to be simply proportional to curvature.

In this paper, we propose a generalisation of the hyperbolic flow model proposed in Refs~\cite{Alias2017,Alias2018} that is applicable to general evolution of tissues that involve fragmentation and fusion processes, and that is independent of spatial dimension. This is achieved by considering two Eulerian fields in space. The first field is a level-set function that represents the interface implicitly as the zero contour level. The second Eulerian field represents the surface density of tissue-synthesising cells at the tissue interface. In a neighbourhood of the tissue interface, the value of this field anticipates influences on cell density that would be exerted along the path to that location. The new coupled partial differential equations governing the evolution of the level set function and of the cell density field provide a new technique to solve hyperbolic curvature flows in complex topological situations. These equations are expressed in a manifestly covariant form independent of space dimension. We note that our method is also applicable to other surface-bound dynamic processes that affect the evolution of an interface, including etching processes~\cite{Adalsteinsson1995,Adalsteinsson1995b}, active membranes~\cite{Ramaswamy2001,Cagnetta2018}, and thin films and foams~\cite{Dexing2009,He2009,Ishida2017,Da2015,Durikovic2001,Zhu2014}.

The paper is organised as follows. In Section~\ref{sec:math-model} we present the mathematical model and derive the governing equations for the level set function and cell density field. Section~\ref{sec:numerical-methods} details the numerical methods employed to solve these equations. In Section~\ref{sec:results-methods}, these methods are compared with analytic solutions and with simulations performed with explicit parameterisations of the interface. We find that a good indicator of numerical accuracy is provided by how well cell density is conserved numerically, i.e., by the discrepancy between expected and numerical total cell numbers. Finally, we apply this new model to situations that could not be modelled by the previous models~\cite{Alias2017,Alias2018}, including tissue fusion, fragmentation, and three dimensions (Section~\ref{sec:results-applications}).

\section{Mathematical model}\label{sec:math-model}
In tissue engineering bioscaffolds, tumour spheroids, epidermeal wound healing, and new bone formation, new tissue is secreted by active cells residing at or near the tissue's surface~\cite{Rumpler2008,Bidan2012,Bidan2013,Jin2018,Vermolen2007,Sutherland1988,Wise2008,Macklin2003,Macklin2005}. We have shown in Refs~\cite{Alias2017,Alias2018} that the evolution of the tissue interface for such surface-localised growth is such that the normal acceleration of the interface depends linearly on curvature. This evolution constitutes a type of hyperbolic curvature flow~\cite{LeFloch2008}.

It is important to emphasise that the presence of this flow is inevitable due to the geometric crowding of new tissue material at concavities of the substrate, and geometric spreading of new tissue material at convexities. Whether actual crowding or spreading occurs during tissue growth depends of course on the presence of other driving forces that may counteract these effects, such as mechanical relaxation, diffusion, and other cell behaviours, but the curvature-induced crowding/spreading force itself is not one that can be turned off. The tissue growth dynamics described by this hyperbolic curvature flow derives directly from the conservation law of the tissue-synthesising cell surface density $\rho$ (number of cells per unit surface) and the fact that the normal velocity of the interface is assumed to be given by
\begin{align}    \label{normal_velocity}
v = \kform \, \rho,
\end{align}
where $k$ is the cell secretory rate (volume of new tissue secreted per cell per unit time)~\cite{Alias2017,Buenzli2015}. In Refs~\cite{Alias2017,Alias2018}, the tissue interface $S(t)$ in two-dimensional space is described by an explicit front-tracking parameterisation $\b\gamma(s,t)$, where $s$ is an arbitrary one-dimensional parameter and $t$ is time. If the parameterisation $\b\gamma(s,t)$ is orthogonal, i.e., time lines $t\mapsto \b\gamma(s,t)$ are perpendicular to the interface everywhere at all times, the equations governing the evolution of the tissue boundary and the cell density are given by~\cite{Alias2017}:
\begin{align}
&\b\gamma_t  = v\b n, \label{gamma_t}
\\&\rho_t
    =
- k\rho^2 \kappa    
+ D \rho_{\ell \ell} 
- A\rho,\label{rho_t}
\end{align}
where subscripts denote partial derivatives, $\ell$ is the arc length such that $\d\ell = g\d s$ with $g=|\b\gamma_s|$, $\b\tau=\b\gamma_s/|\b\gamma_s|$ is the unit tangent vector to $S(t)$, $\b n$ is the outward unit normal to the tissue substrate, and $\kappa = \b\tau\cdot\b n_\ell$ is the signed curvature, such that $\kappa>0$ where the tissue substrate is convex, and $\kappa <0$ where the tissue substrate is concave. In Eq.~\eqref{rho_t}, the term $D\rho_{\ell\ell}$ represents cell diffusion parallel to the interface with diffusivity $D$, and the term $-A\rho$ represents the depletion of active cells at rate $A$.

If the secretory rate $k$ is constant, cell density $\rho$ can be substituted for normal velocity $v$ (and vice versa) by Eq.~\eqref{normal_velocity}, and Eq.~\eqref{rho_t} is equivalent to
\begin{align}
v_t 
    =
- v^2 \kappa    
+ D v_{\ell \ell} 
- Av,\label{v_t}
\end{align}
where $v_t$ corresponds to the normal acceleration of the interface, $\b\gamma_{tt}\cdot\b n$. The first term in the right hand side therefore makes this flow a hyperbolic curvature flow.

The disadvantage of these equations is that they involve an explicit parameterisation of the interface. Explicit parameterisations make it difficult to capture the evolution of tissues undergoing complex changes in morphology. To describe arbitrary interface shapes and changes in topology, we now represent the interface $S(t)$ implicitly as the zero contour of a time-dependent scalar field $\phi$, called the level set function:
\begin{align}
S(t) = \{\b r\ |\ \phi(\b r,t)=0\}.
\end{align} 
The equation that governs the evolution of the level set function $\phi$ is found by differentiating $ \phi(\boldsymbol{\gamma}(s,t),t)=0$ with respect to $t$. Utilising the fact that the unit normal vector of contour levels of $\phi$ is $\boldsymbol{n}=\nabla \phi/|\nabla \phi|$ and that $\boldsymbol{\gamma}_t \cdot \boldsymbol{n}=v$, one gets~~\cite{Sethian1999,Osher2003}
\begin{align}
\phi_t + V|\nabla \phi| = 0,    \label{phi}
\end{align}
where $V$ represents the normal velocities of all the contour lines of $\phi$, and must coincide with the normal velocity $v=k\,\rho$ at the interface. In many applications of the level set method, the normal velocity of the interface is known algebraically, such as in mean curvature flow, where the normal velocity $v$ can extrapolated in a neighbourhood of the interface as $V \propto \nabla\cdot \b n = \nabla\cdot \frac{\nabla \phi}{|\nabla \phi|}$. In these cases, Eq.~\eqref{phi} is the only equation to solve. Even though one is in principle only interested in how the interface evolves, i.e, in how the zero contour level $\phi(\b\gamma(s,t),t)$ evolves, Eq.~\eqref{phi} is solved for $\phi(\b r, t)$ in the whole Cartesian space $\b r$ (or in a restricted band around the interface) using regular PDE techniques~\cite{Sethian1999,Osher2003}, and the zero contour of $\phi$ is determined afterwards.

In our situation, the normal velocity is a solution of a differential equation, Eq.~\eqref{v_t}, which represents the effect of dynamic processes confined to the interface only. Now that the interface is described implicitly, an alternative description of the normal velocity, or cell surface density, that does not refer to the explicit parameter~$s$, is also required. Proceeding similarly to the derivation of the level set equation~\eqref{phi}, we seek a scalar field $\hat\rho(\b r, t)$ that coincides with $\rho(s,t)$ at any point $\boldsymbol{\gamma}(s,t)$ of the interface. Writing
\begin{align}\label{rho}
    \hat\rho\left(\b\gamma(s,t),t\right) = \rho(s,t)
\end{align} 
and differentiating with respect to $t$ gives, after using Eqs~\eqref{gamma_t},~\eqref{rho_t}, and~\eqref{rho},
\begin{align}\label{level-set-like-rho}
    \hat\rho_t + k \hat\rho\b n \cdot \nabla \hat\rho = -k\hat\rho^2 \nabla\cdot\b n + D \nabla_\ups^2 \hat\rho - A \hat\rho,
\end{align} 
where $\nabla_\ups^2\hat\rho=\nabla_\ups\cdot\nabla_\ups \hat\rho$ is the Laplace-Beltrami operator, and where the unit normal $\b n$ and mean curvature $\kappa$ are extended to a neighbourhood of the interface via the level set function~\cite{Sethian1999}:
\begin{align}\label{n-kappa}
\b n=\frac{\nabla\phi}{|\nabla\phi|}, \qquad \kappa=\tfrac{1}{d-1}\nabla\cdot \b n.
\end{align} 
Note that the expression for the mean curvature depends on the spatial dimension $d$ ($d=2$ in Eq.~\eqref{gamma_t})~\cite{Sethian1999}, although some authors define mean curvature as $\nabla\cdot\b n$ irrespective of the dimension, e.g.~\cite{Xu2003}. The sign of $\kappa$ is consistent with our convention provided that $\phi<0$ inside the tissue and $\phi>0$ outside the tissue. The Laplace-Beltrami operator can also be extended to a neighbourhood of the interface. The surface gradient $\nabla_\ups \hat\rho$ can be obtained as the orthogonal projection of $\nabla \hat\rho$ onto the surface, i.e.
\begin{align*}
\nabla_\ups \hat\rho = (\mathbb{I}-\b n\b n^T)\nabla \hat\rho,
\end{align*}
and the surface divergence $\nabla_\ups\cdot \b F$ of a vector field $\b F$ as the trace of the Jacobian matrix $\nabla \b F$ restricted to the tangent plane, or equivalently, as the trace minus the normal component~\cite{Arnoldus2006}:
\begin{align*}
    \nabla_\ups\cdot \b F = \Tr_\ups (\nabla\b F) = \Tr(\nabla\b F) - \b n^T \nabla\b F \b n = \nabla\cdot \b F - \b n^T\nabla\b F\b n.
\end{align*} 
In two-dimensional space ($d=2$), these definitions are consistent with the vectorial representation $\nabla_\ups = \b\tau \frac{\p}{\p\ell}$, where $\frac{\p}{\p\ell} = \b\tau\cdot\nabla$, as expected~\cite{Redzic2001}. With these extensions of the surface gradient and surface divergence away from the interface, the Laplace-Beltrami operator can be calculated as~\cite{Xu2003}
\begin{align}\label{laplace-beltrami}
    \nabla_\ups^2 \hat\rho = \nabla^2 \hat\rho - (\nabla\cdot\b n) (\b n\cdot \nabla \hat\rho) - \b n^T H(\hat\rho) \b n,
\end{align}
where $H(\hat\rho)=\nabla\nabla \hat\rho$ is the Hessian matrix $\{\frac{\p^2\hat\rho}{\p x_i\p x_j}\}$ of $\hat\rho$.

Equation~\eqref{level-set-like-rho} with Eqs~\eqref{n-kappa} and~\eqref{laplace-beltrami} corresponds to an extension of Eq.~\eqref{rho_t} to the whole Cartesian space, that is consistent with the requirement that at the interface, cell density evolves by Eq.~\eqref{rho_t}. The solution $\hat\rho$ at a point $\b r$ outside the interface can be interpreted from Eq.~\eqref{rho} as the anticipated value that the density would take considering the crowding/spreading influence of curvature, the cell diffusion, and the cell depletion effects that would be exerted on a path going from the current interface to the point $\b r$.

An extension $V(\b r, t)$ of the normal velocity $v$ can be defined in the Cartesian space similarly. Setting $V(\b\gamma(s,t), t) = v(s,t)$ and differentiating with respect to $t$ as above, or simply defining
\begin{align}\label{V-hatrho}
V(\b r, t) = k\hat\rho(\b r, t)
\end{align}
based on Eq.~\eqref{normal_velocity}, the partial differential equation (PDE) for $V$ that corresponds to Eq.~\eqref{v_t} is given by
\begin{align}\label{level-set-like-V}
    V_t + V\b n \cdot \nabla V = - V^2\nabla\cdot \b n + D \nabla_\ups^2 V - A V.
\end{align} 

\begin{figure}[t!]
	
	\centerline{
		\includegraphics[trim={0 0 0 0 },width=0.98\linewidth,clip]{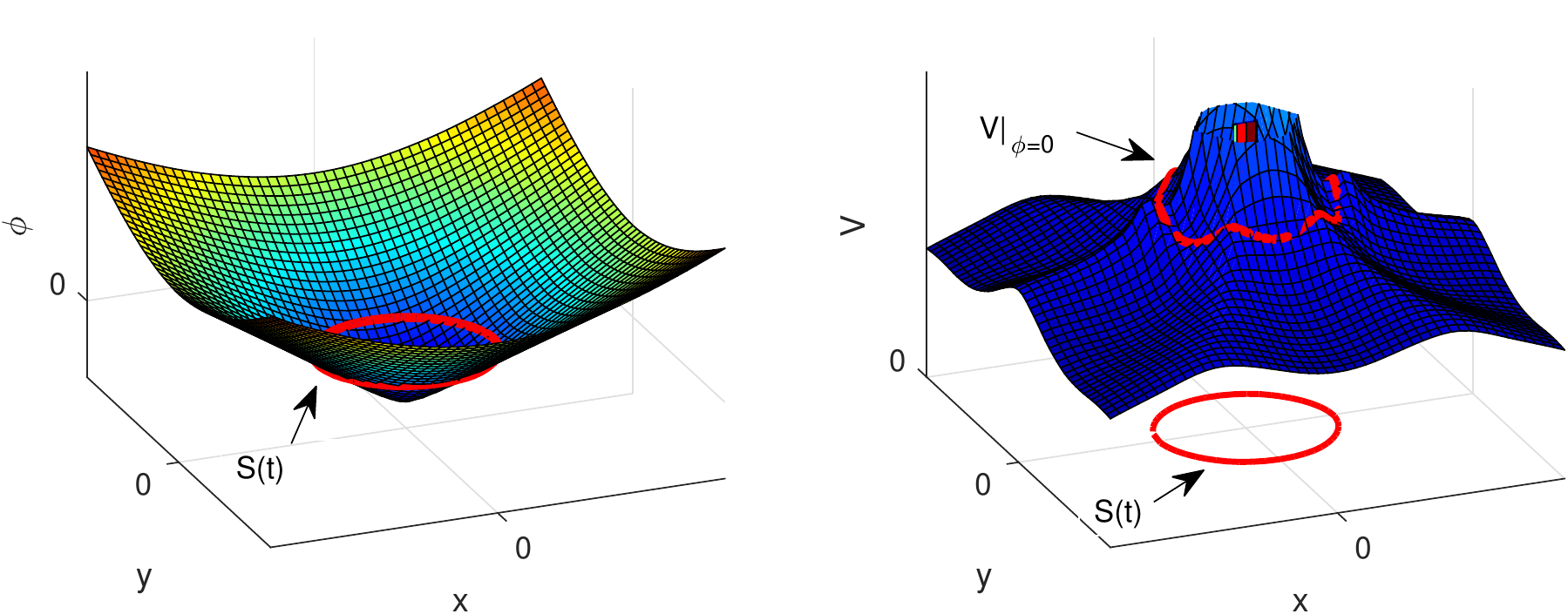}
	}

	\caption{Level set function $\phi$ with interface $S(t)$ at its zero level set (left). The values of velocity function $V$ at the interface points $S(t)$ are the interface velocity $\left. V \right|_{S(t)} =\{V(x,y,t)|(x,y)\in S(t)\}$ (right).}
	\label{fig2}
\end{figure}

Equations~\eqref{phi} and~\eqref{level-set-like-rho} form a system of two nonlinear PDEs that describe implicitly both the position of the tissue interface, and the cell density. For simplicity, in the following, we will always consider $k$ to be a constant so that cell density $\hat\rho$ can be substituted for the velocity field $V$ by Eq.~\eqref{V-hatrho}. The zero contour of $\phi$ provides the set of all points belonging to the interface $S(t)$. Evaluating $V$ at these points then provides the normal velocity of the interface, see Fig.~\ref{fig2}.

While the derivation proposed for Eqs~\eqref{phi} and~\eqref{level-set-like-rho} is based on an explicit parameterisation of the interface in two-dimensional space ($d=2$), the level set equation~\eqref{phi} has the same form in higher dimensions~\cite{Sethian1999}, and Eq.~\eqref{level-set-like-rho} matches the evolution equation of a surfactant $\Gamma$ on a moving boundary in two or three dimensions ($d=2,3$)~\cite{Wong1996,Xu2003,James2004}. Equations~\eqref{phi} and~\eqref{level-set-like-rho} thus generalise the hyperbolic curvature flow model of tissue growth of Ref.~\cite{Alias2017} to three-dimensional space. Since these equations are expressed in terms of vector calculus operators, these equations are also manifestly covariant with respect to changes of reference frames. In contrast to the passive evolution of surfactants on moving boundaries~\cite{Xu2003}, here surface cell density plays an active role in the evolution of the interface: the evolution of the density influences the evolution of the interface by Eq.~\eqref{V-hatrho}. This results in a strongly coupled system which reflects the mechanistic curvature control of tissue growth and the hyperbolic character of the flow. Coupling occurs via $V$ in Eq.~\eqref{phi}, and via $\b n$ and $\kappa$ (expressed in terms of $\phi$) in Eq.~\eqref{level-set-like-V}.

\section{Numerical methods}\label{sec:numerical-methods}


Having derived evolution equations that allow us to represent the tissue interface and the cell density in an implicit manner robust to topological changes of the interface, we now devise several possible strategies to solve these equations numerically, drawing on existing numerical methods developed for the level-set method and for PDEs on moving boundaries~\cite{Sethian1999,Osher2003,Peng1999}.

In standard situations where the velocity field is known algebraically, the level set method is known to be more accurate when the level set function is re-initialised as a signed distance function with $|\nabla\phi|=1$~\cite{Sethian1999}. It is also known that maintaining the signed distance property of $\phi$ can be achieved by extrapolating the velocity field in the orthogonal direction, which is achieved by solving
\begin{align}\label{V-ortho-ext}
\nabla V \cdot \nabla \phi = 0
\end{align} 
for $V$~\cite{Sethian1999,Osher2003}. In contrast, the velocity field $V$ solution of the differential equation~\eqref{level-set-like-V} anticipates values at future locations of the interface by accounting for the curvature-induced acceleration or deceleration of the interface. The different numerical strategies we devise differ by (i) whether or not the level set function~$\phi$ is re-initialised as a signed distance function; and (ii) whether or not the velocity field~$V$ is re-initialised by the orthogonal extrapolation~\eqref{V-ortho-ext} after determining its value at the interface using Eq.~\eqref{level-set-like-V}:
\begin{description}
\item[Method 1] No re-initialisation of $\phi$ nor of $V$;
\item[Method 2] Re-initialisation of $\phi$ as a signed distance function; no re-initialisation of $V$;
\item[Method 3] Re-initialisation of $\phi$ as a signed distance function; re-initialisation of $V$ by orthogonal extrapolation.
\end{description}
Methods~1--3 therefore explore a trade-off between using a possibly more accurate determination of velocity (no re-initialisation of $V$, Methods~1,2), and using a signed distance function for the level set function (re-initialisation of $\phi$, Methods~2,3). Methods~1 and~3 are consistent, in the sense that all the contour levels of $\phi$ are evolved by the velocity field $V$. Method~2 is not consistent in this sense, because  while all the contour levels of $\phi$ are initially evolved by the solution $V$ to Eq.~\eqref{level-set-like-V}, these contour levels are subsequently reorganised to re-initialise $\phi$ as a distance function. However, Method~2 may potentially combine the advantages of dealing with a signed distance function for $\phi$, and anticipating future values of $V$ away from the current interface.

The general solution algorithm we use to solve Eqs~\eqref{phi} and~\eqref{level-set-like-V} jointly is based on operator splitting and discrete time stepping, and is summarised as follows:
\begin{description}
    \item[Step 1] \emph{Initialisation.} $\phi$ is set as a signed distance function $\phi^0$ to the initial interface $S(0)$, and $V$ is initialised to be a uniform constant $v^0$ along the interface. In Methods 1 and 2, the velocity field $V$ is extended away from the interface by solving Eq.~\eqref{level-set-like-V} with $\phi^0$ fixed until convergence (see Step 3). In Method 3, the orthogonal extrapolation of the initial interface velocity is achieved by simply setting $V$ to $v^0$ in the whole computational domain.
    \item[Step 2] \emph{Level-set function update.} The level set function at time step $n$, $\phi^n$, is evolved using explicit stepping in time, leading first to a temporary update $\phi^{n+1/2}$. In general, $\phi^{n+1/2}$ is no longer a signed distance function. Re-initialisation of $\phi^{n+1/2}$ to a signed distance function is performed in Methods~2 and~3. This leads to the full time step update $\phi^{n+1}$.
    \item[Step 3] \emph{Velocity field update.} The velocity field at time step $n$, $V^n$, is evolved by solving Eq.~\eqref{level-set-like-V} using a semi-implicit time stepping scheme, leading first to a temporary update $V^{n+1/2}$. In Methods~1--2, $V^{n+1/2}$ corresponds to the full time step update, i.e. $V^{n+1}=V^{n+1/2}$. In Method~3, $V^{n+1/2}$ is re-initialised by extrapolating its value at the interface $\phi^{n+1}=0$ in the orthogonal direction by solving Eq.~\eqref{V-ortho-ext}, leading to the full time step update $V^{n+1}$.
\end{description}
The zero level set $\phi^{n+1}=0$ provides the new location of the interface $S^{n+1}$, and the new interface velocity is provided by evaluating the field $V^{n+1}$ at this location. Steps~2 and~3 are repeated iteratively to evolve the solutions $\phi$ and $V$ in time.

We now describe the numerical discretisation algorithms involved in these steps in more detail. These algorithms are based on Ref.~\cite{Sussman1994} with some modifications made to account for the discretisation of the Laplace--Beltrami operator. We restrict the formulas to two dimensions for simplicity.

\paragraph{Discretisation of gradients}
Equations~\eqref{phi} and Eq.~\eqref{level-set-like-V} involve the spatial gradient operator $\nabla$, which we discretise using upwinding based on the velocities $V n_1$ along $x$ and $V n_2$ along $y$, where $n_1$ and $n_2$ are the components of the unit normal $\b n = (n_1, n_2)$~\cite{Sethian1999,Osher2003}. I.e., $\nabla \phi = (\phi_x, \phi_y)$, with
\begin{align}\label{grad}
    \phi_x = 
    \begin{cases}
        \phi_x^-, & \text{if } V n_1 > 0
        \\\phi_x^+, & \text{if } V n_1 \leq 0
    \end{cases},
    &\qquad \phi_y = 
    \begin{cases}
        \phi_y^-, & \text{if } V n_2 > 0
        \\\phi_y^+, & \text{if } V n_2 \leq 0
    \end{cases},
\end{align} 
where $\phi_x^-$, $\phi_x^+$, $\phi_y^-$, and $\phi_y^+$ are backward ($-$) and forward ($+$) high-resolution Hamilton--Jacobi weighted essentially non-oscillatory (HJ-WENO) discretisations of the partial derivatives (and likewise for $\nabla V$)~\cite{Jiang2000,Osher2003}. These discretisations are fifth order accurate in smooth regions of $\phi$ and $V$ but revert to lower order when interpolating across singularities, which occur for example after the emergence of cusps in the interface~\cite{Alias2017}. The term $V|\nabla\phi|$ in Eq.~\eqref{phi} involves the norm of the gradient, and is discretised using Godunov's method with HJ-WENO discretisations~\cite{Peng1999,Osher2003}.

\paragraph{Normal vector and curvature}
At cusps of $\phi$, the unit normal vector and curvature in Eq.~\eqref{n-kappa} are ill-defined. To alleviate the problem of the discontinuity of $\b n$ for the numerical scheme, we follow Ref.~\cite{Sethian1999} and define the unit normal vector $\b n = (n_1, n_2)$ by normalising the average of the four limiting normal vectors that can be calculated by the HJ-WENO backward and forward discretisations:
\begin{align*}
    \b n_{1,2,3,4} = \frac{(\phi_x^\pm, \phi_y^\pm)}{\sqrt{(\phi_x^\pm)^2 + (\phi_y^\pm)^2}}.
\end{align*}

In contrast, we use second order central finite difference for $\phi_x$, $\phi_y$, $\phi_{xx}$, $\phi_{yy}$, and $\phi_{xy}$ in the numerical evaluation of curvature in the formula~\cite{Sethian1999}
\begin{align}
\kappa = \text{cl}\left(\frac{1}{d-1}\frac{\phi_y^2\phi_{xx} - 2\phi_x \phi_y \phi_{xy} + \phi_x^2\phi_{yy}}{(\phi_x^2+\phi_y^2)^{3/2}}\right),
\end{align} 
where 
$\text{cl}(\xi)=
\max({\kappa_\text{min}, \min(\xi, \kappa_\text{max})})
$ is a clamping function that enforces the computed signed curvature $\kappa$ to remain between the minimum value $\kappa_\text{min} = -1/\Deltaup x$ and the maximum value $\kappa_\text{max} = 1/\Deltaup x$ for a spatial discretisation step $\Deltaup x$ \cite{Osher2003}. 

\paragraph{Level-set function time stepping}
We solve the level-set equation~\eqref{phi} numerically with a simple first-order forward Euler discretisation in time for convenience. Level set methods are known to be more sensitive to spatial accuracy than temporal accuracy~\cite{Osher2003}. We have found that using a third order total variation diminishing
Runge--Kutta method \cite{Shu1988} adds significantly more computation time and complexity without changing the results significantly. Our aim is to compare Methods 1--3 irrespective of the time discretisation scheme, so that below we only report results from the simpler first-order Euler discretisation. With this explicit time stepping scheme, the update from $\phi^n$ to $\phi^{n+1/2}$ is given by
\begin{align*}
    \phi^{n+1/2} = \phi^n - \Deltaup t V^n|\nabla \phi^n|,
\end{align*}
with current values $\b n^n$ and $V^n$ of the unit normal and velocity field, and a time increment $\Deltaup t$~\cite{Sethian1999,Osher2003}.

\paragraph{Level-set function re-initialisation}
The level-set function update $\phi^{n+1/2}$ may no longer represent the signed distance function to the interface, even though the zero contour $\phi^{n+1/2}=0$ represents the new location of the interface. In Methods~2 and~3, the level-set function $\phi^{n+1/2}$ is re-initialised to a signed distance function $\phi^{n+1}$ by iterating
\begin{align}\label{reinit}
    \psi^{\nu+1} = \psi^\nu - \Deltaup t\  \mathbb{S}\,(\psi^\nu)\left(|\nabla \psi^\nu| - 1\right)
\end{align} 
over $\nu=0,1,\ldots$ to steady state $\psi^\infty=\phi^{n+1}$, starting from the initial condition $\psi^0 = \phi^{n+1/2}$, where
\begin{align}   \label{sign_function}
\mathbb{S}(\psi) = \frac{\psi}{\sqrt{\psi^2 + |\nabla \psi|^2 (\Deltaup x)^2}}
\end{align} 
is a smoothed sign function~\cite{Peng1999,Sussman1994}. This iterative approach corresponds to finding the steady state of $\psi_\tau = - \mathbb{S}(\psi)(|\nabla \psi|-1)$ with respect to the virtual time $\tau$, which occurs when $\psi$ is a signed distance function. Because $\mathbb{S}(0)=0$, the zero contour is unaffected by this re-initialisation procedure. We use the stopping criterion $\frac{1}{M}\sum_{(i,j): |\psi_{ij}^\nu|<\beta} \left| |\nabla\psi_{ij}^\nu|-1\right| \leq \epsilon_\text{reinit}\Deltaup x\Deltaup y$, where $\beta$ is the one-sided width of the band around the interface within which we require the signed distance function to be accurate, and $M$ is the number of spatial discretisation points $(i,j)$ in the sum. In practice, we choose $\beta=5\Deltaup x$ except for Step 1 (initialisation), where we choose $\beta=20\Deltaup x$.

\paragraph{Velocity field time stepping}
To solve Eq.~\eqref{level-set-like-V}, we rewrite it in the form $V_t = D\nabla^2 V + \alpha$, where $\nabla^2V$ is the isotropic diffusion contribution of $\nabla_\ups^2V$ in Eq.~\eqref{laplace-beltrami}, and
\begin{align}\label{alpha}
  \alpha = &- V \b n \cdot \nabla V - (d-1) \kappa V^2 - D(d-1)\kappa \b n \cdot \b \nabla V \notag
\\&- D \left(n_1^2 V_{xx} + 2 n_1 n_2 V_{xy} + n_2^2 V_{yy}\right) - A V
\end{align}
contains all the other contributions in Eq.~\eqref{level-set-like-V}. We use a simple semi-implicit scheme in which $D\nabla^2V$ is solved implicitly using the alternative direction implicit method (ADI) \cite{Press1992}, and the terms in $\alpha$ are solved explicitly using first-order forward Euler discretisation as above. The ADI is a two-step method, so that the velocity time step update $V^n\to V^{n+1/2}$ is given by solving sequentially
\begin{align}
  &\frac{\widetilde{V}^{n+1/2} - V^n}{\Deltaup t/2} = D \left( \widetilde{V}_{xx}^{n+1/2} + V_{yy}^{n} \right) + \alpha^n,\label{predictor}
  \\&\frac{V^{n+1/2} - \widetilde{V}^{n+1/2}}{\Deltaup t/2} = D \left( \widetilde{V}_{xx}^{n+1/2} + V_{yy}^{n+1/2} \right) + \alpha^n.\label{corrector}
\end{align} 
Equation~\eqref{predictor} is the predictor step determining $\widetilde{V}^{n+1/2}$ implicitly, and equation~\eqref{corrector} is the corrector step determining $V^{n+1/2}$ implicitly. The derivatives $V_{xx}$, $V_{yy}$, $V_{xy}$, and $\widetilde{V}_{xx}$ in Eqs~\eqref{predictor},~\eqref{corrector}, and~\eqref{alpha} are discretised with second-order central difference, whereas first-order derivatives in Eq.~\eqref{alpha} use upwind HJ-WENO discretisations. Each of the predictor and corrector steps requires solving a tridiagonal matrix system with constant matrix coefficients.

\paragraph{Velocity field re-initialisation}
In Method~3, the velocity field $V^{n+1/2}$ is re-initialised to a velocity field $V^{n+1}$ by extrapolating the interface values of $V^{n+1/2}$ in the orthogonal direction, so that it satisfies Eq.~\eqref{V-ortho-ext}. This is achieved by iterating
\begin{align}
    W^{\nu+1} = W^\nu - \Deltaup t \mathbb{S}(\phi^n) \b n^n\cdot \nabla W^\nu
\end{align}
over $\nu=0,1,\ldots$ to steady state $W^\infty = V^{n+1}$, starting from the initial condition $W^0 = V^{n+1/2}$~\cite{Osher2003}. In practice, it is sufficient to perform 10 iterations every 10$^\text{th}$ time step for $V$ to be well extrapolated orthogonally in the band $\beta=5\Delta x$ each side of the interface.

\paragraph{Simulation parameters}
All the simulations are performed with an initial interface velocity of $v^0=0.016\,\text{mm/day}$~\cite{Alias2017}, except simulations of trabecular bone formation and resorption, which use $v^0=0.001\,\text{mm/day}=1\,\muup\text{m/day}$~\cite{Martin1998}. The computational domain is chosen to extend at least 50\% more than the maximum diameter of the interface in all cardinal directions. Space discretisation and the tolerance parameter for reinitialisation are chosen depending on the size of the initial pore. Because of the explicit time stepping scheme, space and time discretisation are restricted by the Courant--Friedrichs--Lewy condition~\cite{Osher2003}. For each figure, these parameters were chosen such that reducing them further would not change the results noticeably. The values of discretisation parameters are listed in the figure captions. Identical discretisation parameters are chosen to compare Methods 1--3.

\section{Results} \label{Section_Results}
We first investigate the numerical accuracy of Methods~1--3 in Sec.~\ref{sec:results-methods}. Application examples of complex geometries, fusion and fragmentation are presented in Sec.~\ref{sec:results-applications}.

\begin{figure*}[t]
	\centering{\includegraphics{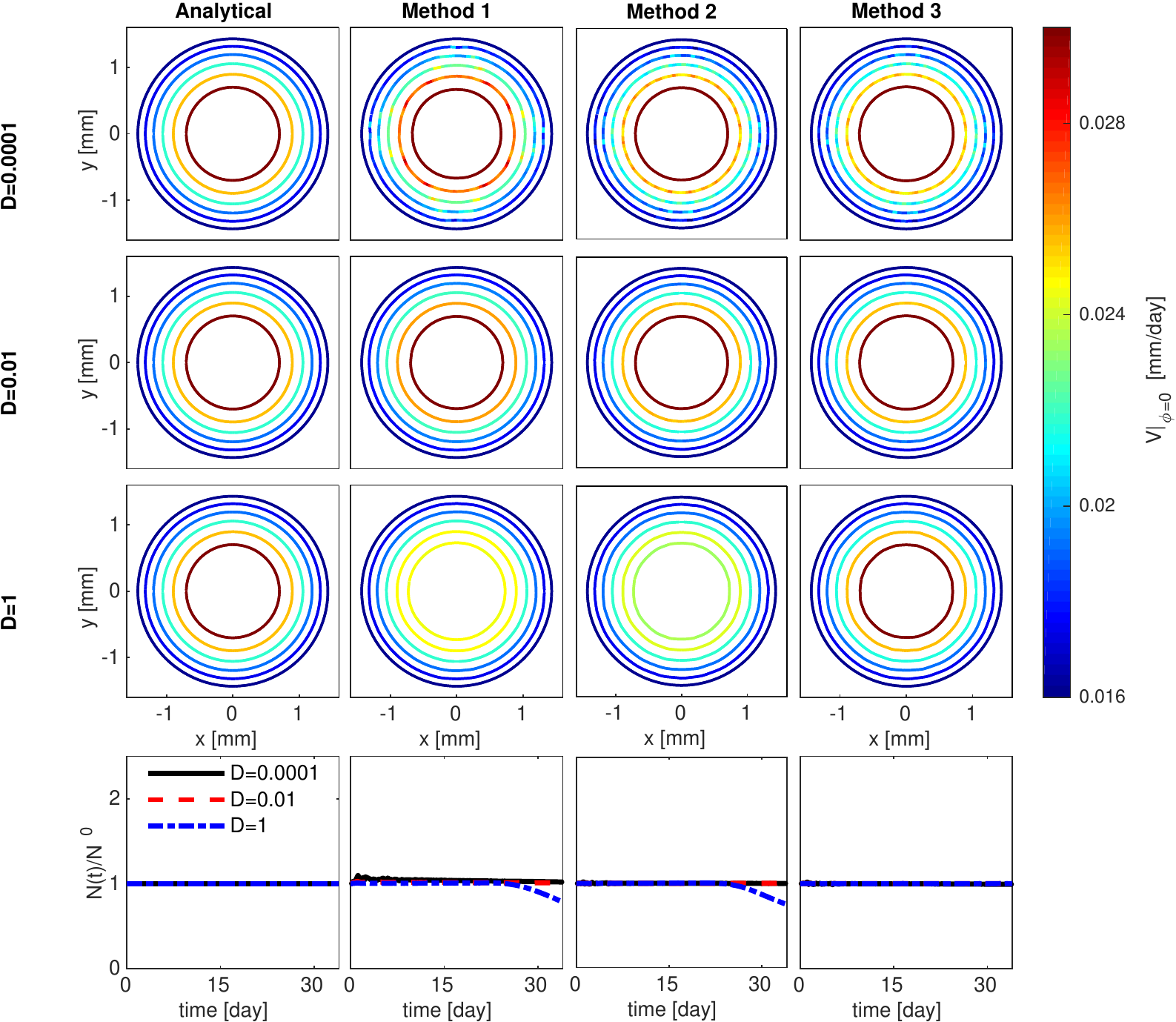}}
	
	\caption{Tissue deposition within a circular pore obtained analytically (first column), and simulated by Methods 1--3 (columns 2--4) with different lateral diffusivities $D$ in $\text{mm}^2/\text{day}$ (rows). The tissue interface is shown at regular time intervals of 6.8 days until 34~days and coloured according to the interface velocity. The evolution of normalised cell number is shown in the last row. Simulation parameters: $\epsilon_\text{reinit}=5$, $\Deltaup x=\Deltaup y=0.0357\,\text{mm}$ and $\Deltaup t=0.017\,\text{days}$.}
	\label{fig:circular}
\end{figure*}

\subsection{Numerical simulations using Methods 1--3}\label{sec:results-methods}
To estimate the accuracy of the numerical methods 1--3, we compare the shape of the interface and the normal velocity at regular time intervals either with analytical expressions (rotation-symmetric solution), or with results obtained using explicit parameterisations of the interface in simple two-dimensional geometries~\cite{Alias2017,Alias2018}. We also check for conservation properties of the tissue-synthesising cells. For simplicity, we assume in this section that the tissue-synthesising cells are not depleted ($A=0$), meaning that total cell number is constant. Total cell number $N(t)$ is estimated numerically by interpolating $V$ at the interface location $\phi=0$, and integrating the interpolation numerically. Normalising by the initial cell number $N^0$ gives:
\begin{align}
    \frac{N(t)}{N^0} = \frac{\int_{S(t)}\hat\rho\,\d\ell}{\int_{S(0)}\hat\rho\,\d\ell} = \frac{1}{v^0S^0}\int_{S(t)}V \d\ell,
\end{align}
where $S^0$ is the initial interface perimeter.

\begin{figure*}[t]
    \centering{\includegraphics{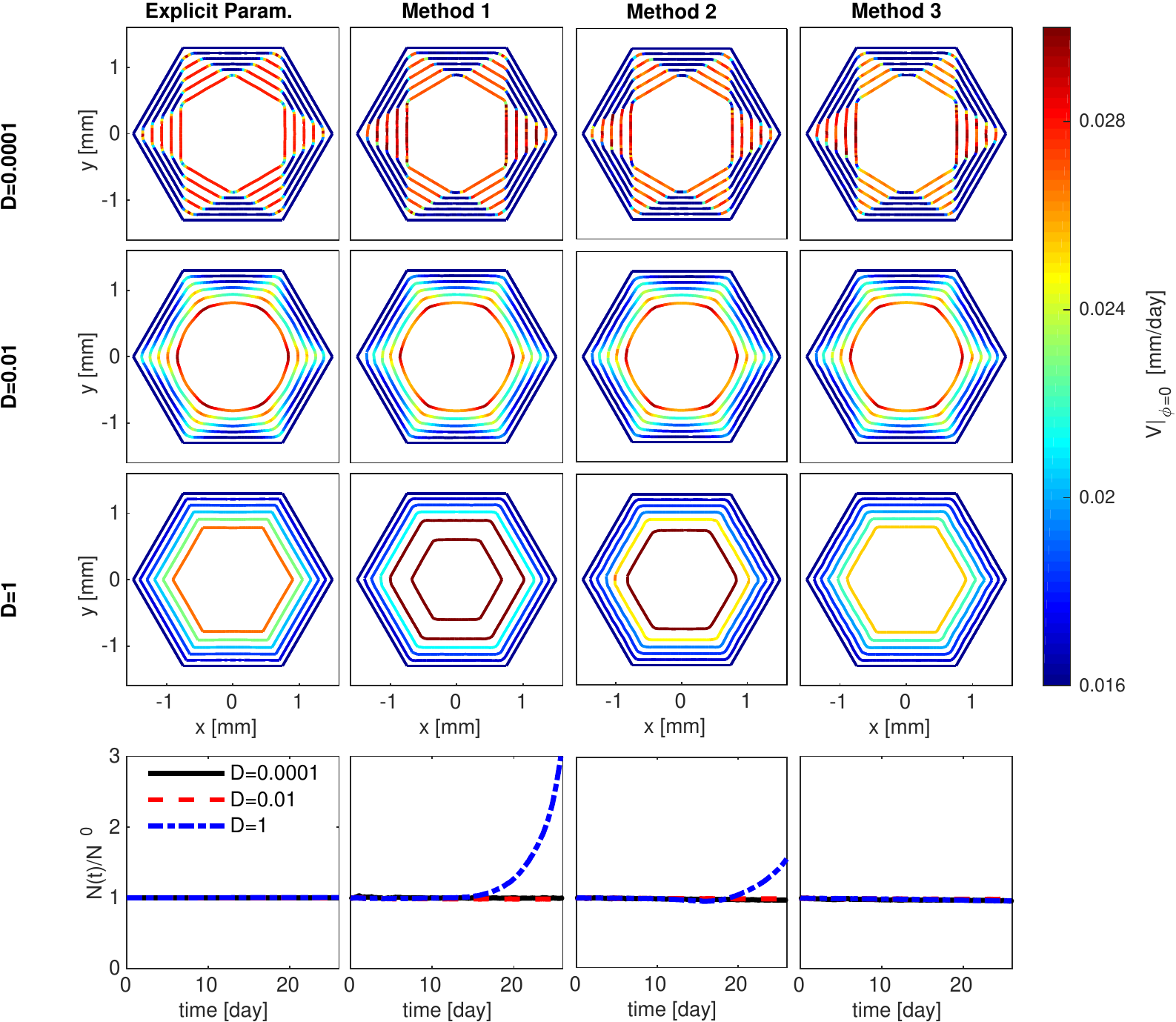}}
	
	\caption{Tissue deposition within a hexagonal pore simulated using an explicit parameterisation (first column)~\cite{Alias2017}, and using Methods 1--3 (columns 2--4) with different lateral diffusivities $D$ in $\text{mm}^2/\text{day}$ (rows). The tissue interface is shown at regular time intervals of 5.2 days until 26~days and coloured according to the interface velocity. The evolution of normalised cell number is shown in the last row. Simulation parameters: $\epsilon_\text{reinit}=5$, $\Deltaup x=\Deltaup y=0.0357\,\text{mm}$ and $\Deltaup t=0.013\,\text{days}$. Explicit parameterisation results are obtained using $\Deltaup \theta=0.0349$ and $\Deltaup t=0.0163\,\text{days}$~\cite{Alias2017}.}\label{fig:hexagonal}
\end{figure*}

\begin{figure*}[t]
    \centering{\includegraphics{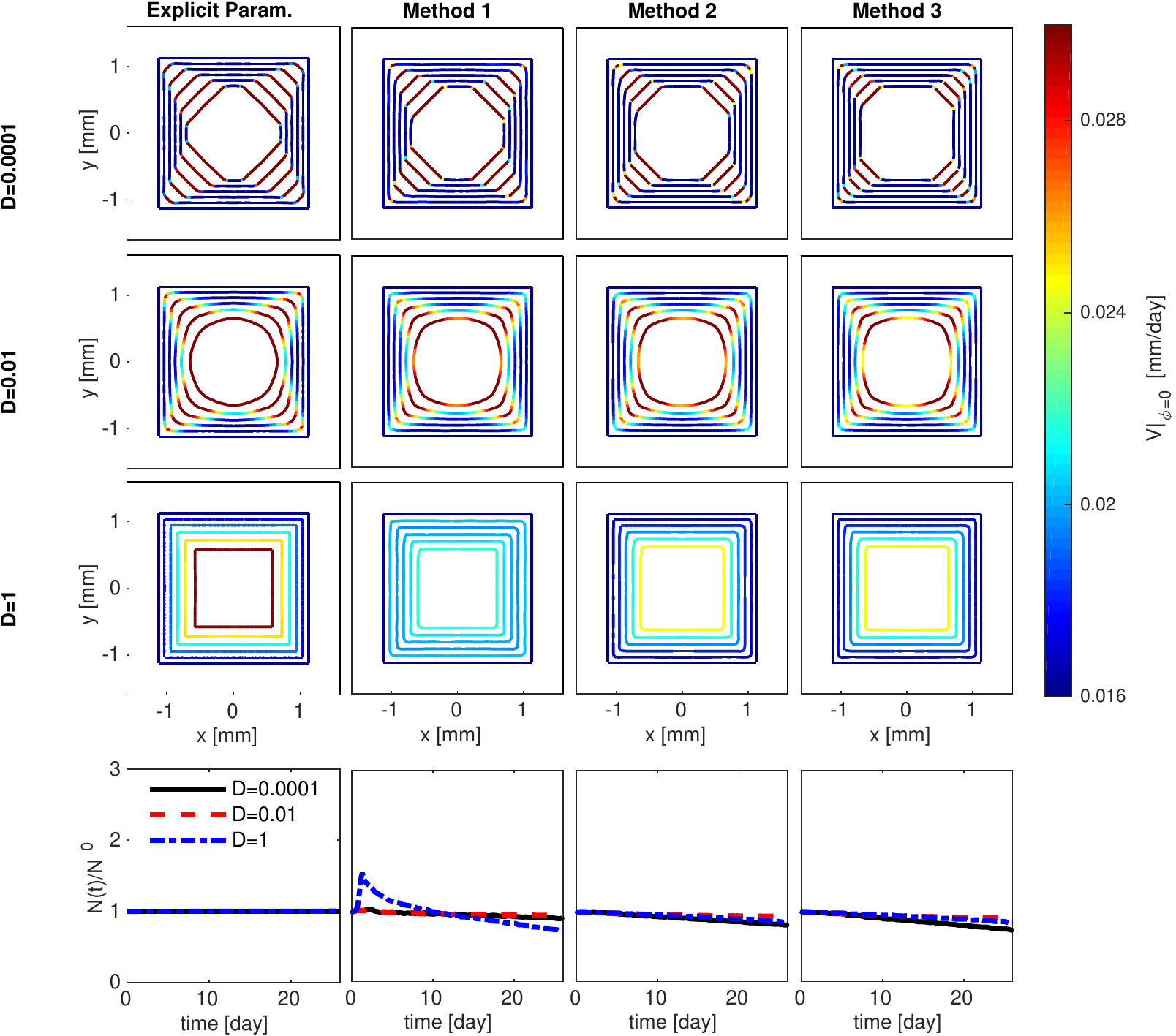}}
	
	\caption{Tissue deposition within a square pore simulated using an explicit parameterisation (first column)~\cite{Alias2017}, and using Methods 1--3 (columns 2--4) with different lateral diffusivities $D$ in $\text{mm}^2/\text{day}$ (rows). The tissue interface is shown at regular time intervals of 5.2 days until 26~days and coloured according to the interface velocity. The evolution of normalised cell number is shown in the last row. Simulation parameters: $\epsilon_\text{reinit}=5$, $\Deltaup x=\Deltaup y=0.0357\,\text{mm}$ and $\Deltaup t=0.013\,\text{days}$. Explicit parameterisation results are obtained using $\Deltaup \theta=0.0196$ and $\Deltaup t=0.0163\,\text{days}$~\cite{Alias2017}.}
	\label{fig:square}
\end{figure*}

We start by investigating the infilling of a circular pore~\cite{Alias2017,Alias2018} (Figure~\ref{fig:circular}), and then consider the infilling of hexagonal and square pores (Figures~\ref{fig:hexagonal}--\ref{fig:square}) to see how Methods 1--3 handle cusps in the interface of increasing sharpness. More realistic geometries are presented in Section~\ref{sec:results-applications}.

\paragraph{Circular pore infilling}
The evolution of the radius $R(t)$ and interface velocity $V(t)=-R_t(t)$ of an infilling circular pore are given by
\begin{align}   \label{circle_analytical}
R(t) = R^0\sqrt{1-2\frac{v^0}{R^0}t}   \quad , \quad
V(t) = v^0\frac{R^0}{R(t)},
\end{align} 
where $R^0$ is the initial radius~\cite{Alias2017}. In the simulations, we choose $R^0=9\,\text{mm}/(2\pi)$ such that the initial pore perimeter is $S^0=9\,\text{mm}$, and we set $v^0 = 0.016\,\text{mm}/\text{day}$ as in Ref.~\cite{Alias2017}. 

By symmetry, cell diffusivity $D$ has no effect on the evolution of the interface and cell surface density in such a situation. However, it is clear from Figure~\ref{fig:circular} that the accuracy of the numerical simulations depends on diffusivity. The simulations are less accurate at low and high diffusivities, to different degrees depending on the method. The accuracy of Method~3 appears to be less sensitive to the degree of diffusivity and matches the analytical result very well. Methods~1 and~2 compare well with the analytical result so long as the total cell number is conserved (Fig.~\ref{fig:circular}, bottom row). There is a significant numerical loss of cells developing at late times in Methods~1 and~2 at high diffusivities that results in an inaccurate evolution of interface and cell density at these times.

\paragraph{Hexagonal and square pore infilling}
Figures~\ref{fig:hexagonal} and~\ref{fig:square} show the infilling of hexagonal and square pores, respectively. The initial pore perimeter is $S^0=9\,\text{mm}$ like in the circular pore case, so that the initial number of tissue-synthesising cells is the same. The first column of Figs~\ref{fig:hexagonal} and~\ref{fig:square} represents simulations obtained by the high-resolution conservative numerical schemes of~\cite{Alias2017}. Like in the circular pore case, all the methods perform well at intermediate diffusivity $D=0.01\,\text{mm}^2/\text{day}$. At low diffusivity $D=0.0001\,\text{mm}^2/\text{day}$, Method~1 performs better than Methods~2 and~3, particularly in the square pore case where corners of the interface are more acute. However, at high diffusivity, Method~1 violates cell conservation significantly, and the interface and velocities obtained by Method~3 are closest to the simulations of Ref.~\cite{Alias2017}.

The numerical violation of cell conservation in Methods~1--3 is more severe initially for the interfaces with more acute cusps. At high diffusivities, the angle of the cusps remains the same throughout the simulation, and this is why cell conservation is harder to achieve numerically. At low diffusivity, cusp angle is initially reduced by a factor two due to the sideways propagation of shock waves~\cite{Alias2017,Alias2018}. Simulations with intermediate diffusivity tend to smooth the interface, which helps conserve cell numbers numerically.

\subsection{Application to complex geometries, fusion and fragmentation}\label{sec:results-applications}
We now consider applications of our model to complex geometries that the conservative Monge parameterisation methods developed in Refs~\cite{Alias2017,Alias2018} cannot represent. The comparison of Methods~1--3 in Section~\ref{sec:results-methods} reveals that numerical conservation can be difficult to achieve, particularly as cusps in the interface emerge. None of Methods~1--3 is explicitly conservative, so that tracking conservation is an important indicator of accuracy. Method~3 resulted in the best conservation properties in general, although Method~1 may still perform better at low diffusivities. Most often, curvature flow models smooth the interface, but the hyperbolic curvature flow considered here only smoothes the interface provided lateral cell diffusion is balanced by curvature-induced crowding or spreading of cells at concavities or convexities of the interface. This depends on the initial interface and diffusivity. 

\paragraph{Bone formation on a single trabecular spicule}
\begin{figure}[tb] \centering{\includegraphics[trim={50 -30 60 0}, width=0.7\columnwidth] {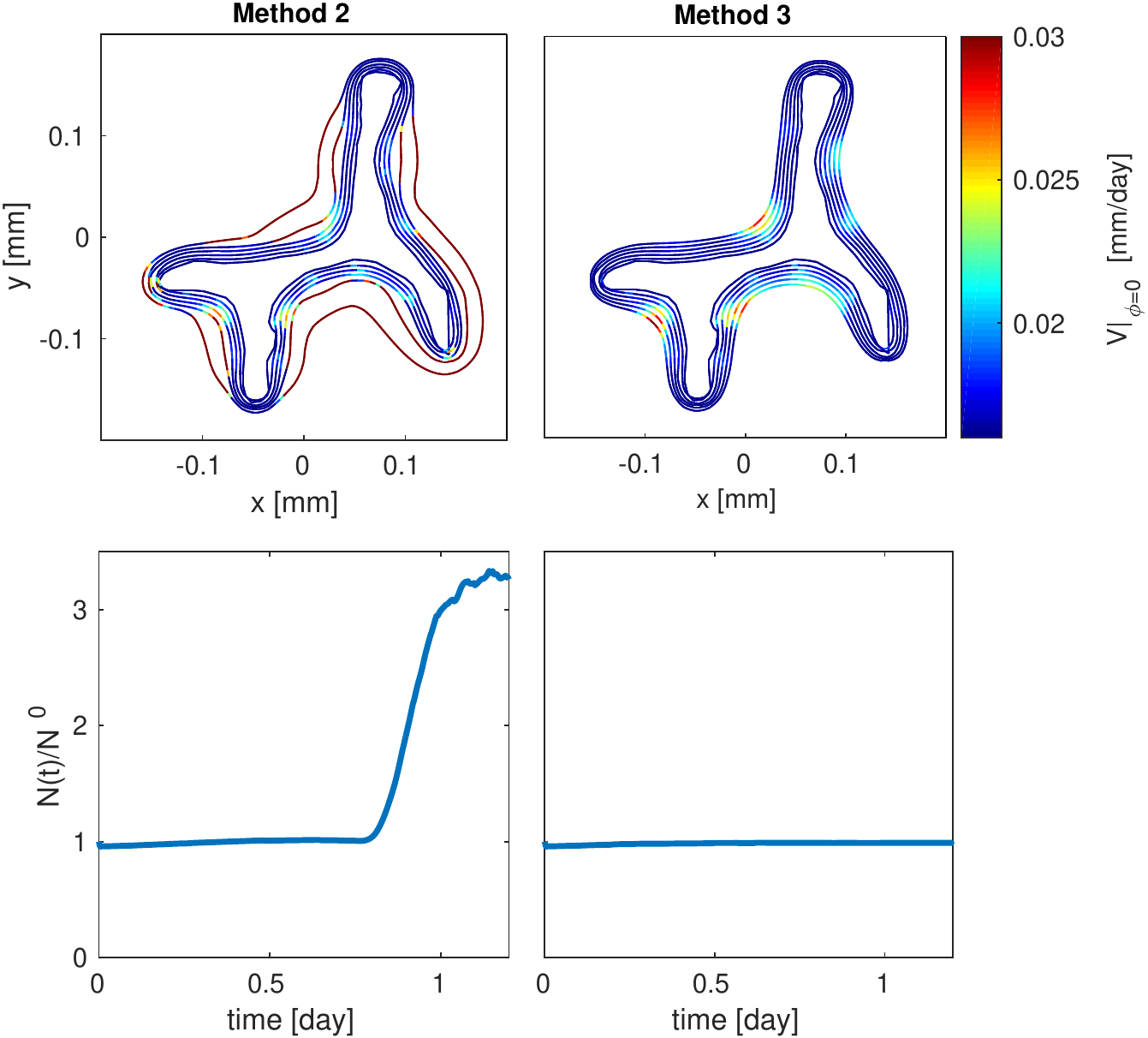}}
    \caption{Bone deposition around a single trabecular spicule obtained using Method~2 (left column) and Method~3 (right column). The initial shape of trabecular spicule is extracted from~\cite[Fig. 2]{Shiga2016}. The interface is shown at regular time intervals of 0.24 days and coloured according to the interface velocity (top row). The evolution of normalised cell number is shown in the bottom row. Simulation parameters: $D=0.0001\,\text{mm}^2/\text{day}$, $\epsilon_\text{reinit}=1000$, $\Deltaup x= \Deltaup y=0.0075\,\text{mm}$, $\Deltaup t=0.00017\,\text{days}$.}
\label{fig:single-trabecula}
\end{figure}
To represent a realistic geometric situation, we take for initial interface the surface of a single trabecular bone spicule seen in an experimental cross section~\cite{Shiga2016}. We then consider the apposition of new bone layers on this surface, as would occur for example by bone mechanical adaptation~\cite{Frost1987,Lerebours2016b}. Figure~\ref{fig:single-trabecula} shows the result of these simulations, where we assume $D=0.0001$ mm$^2/$day and $A=0$. Since Figures~\ref{fig:circular}--\ref{fig:square} show that cell conservation serves as an important indicator of accuracy of interface motion and interface velocity, and since there are no conservative simulations to compare with in this geometry, we now use cell conservation as an indicator of accuracy to compare Methods~1--3.

We see from Figure~\ref{fig:single-trabecula} (bottom row) that Method~2 violates the conservation of cells significantly, greatly overestimating the amount of new bone layers produced. Similar results hold for Method 1, in which cell conservation is also violated significantly (not shown). In contrast, Method~3 maintains cell numbers within 98.1\% of their initial value.

In the remainder of the paper, we use Method~3 to illustrate the capabilities of the level-set formulation of the hyperbolic curvature flow model of tissue growth. We use this method to model complex evolving geometries, including fusion and fragmentation of tissues, and we rate numerical accuracy by tracking cell density conservation properties.

\paragraph{Fusion of two trabecular bone struts and time irreversibility}
\begin{figure} \centering\includegraphics[trim={0 0 0 0}, width=1\columnwidth]{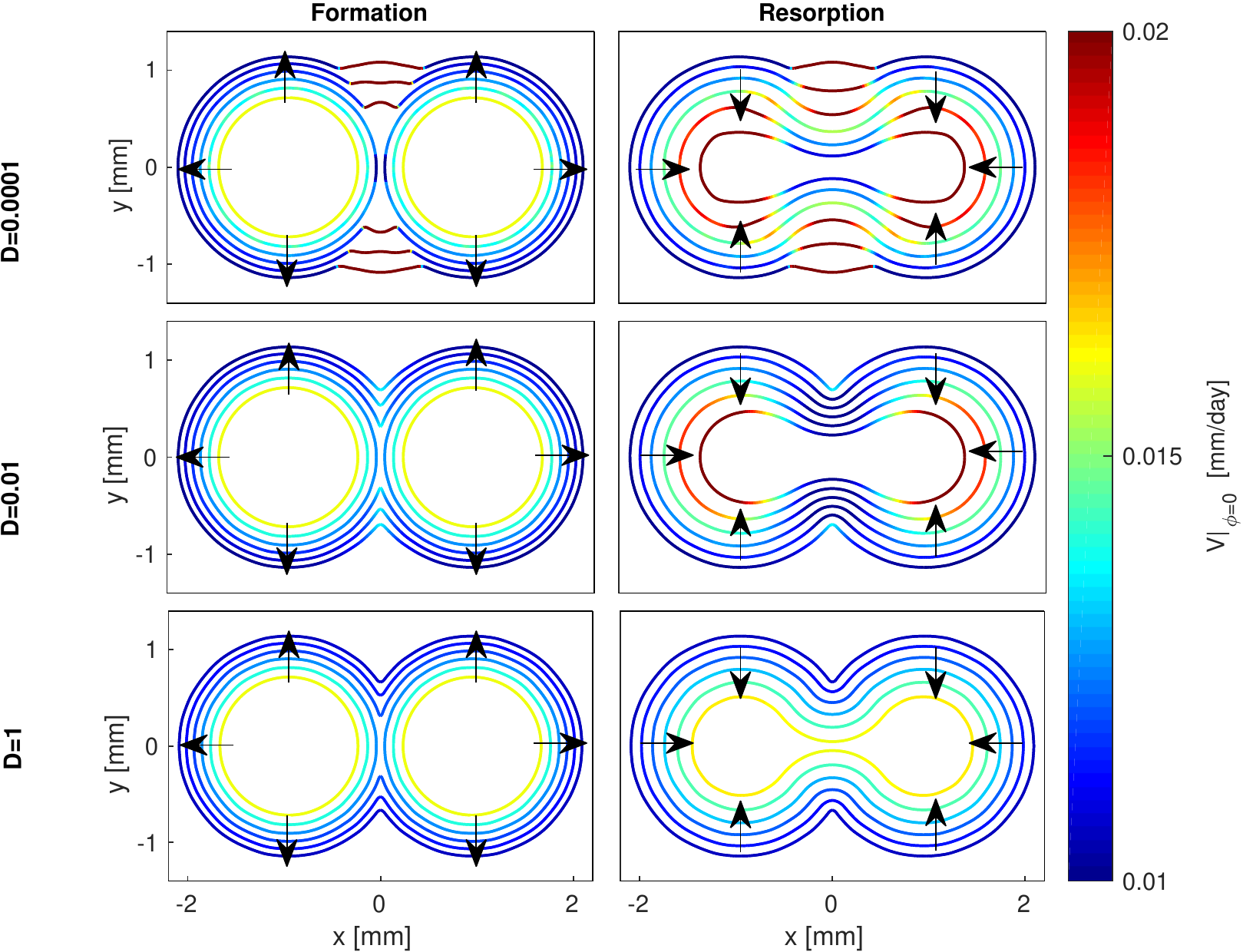}
\caption{Evolution of two trabecular struts (seen in cross-section) during outward tissue deposition (left column), and subsequent inward tissue resorption (right column) with different diffusivities. The initial state for outward motion is two separate trabecular bone struts with a normal velocity of 0.016\,mm/day (yellow-green), which start to merge at $t=17\,\text{days}$. The initial state for resorption is the final state of tissue deposition. Arrows indicate time evolution. Simulation parameters: Method~3, $A=0$, $\epsilon_\text{reinit}=10$; Formation: $\Deltaup x=\Deltaup y=0.0278\,\text{mm}$, $\Deltaup t=0.0085\,\text{day}$; Resorption: $\Deltaup x=\Deltaup y=0.0312\,\text{mm}$, $\Delta t=0.0123\,\text{days}$.}
\label{fig_circle_merge}
\end{figure}
To model a situation where there is a topological change in the interface, we consider the fusion of two trabecular bone struts seen in cross-section (Figure~\ref{fig_circle_merge}). As in Section~\ref{sec:results-methods}, we take the total perimeter to be $9\,\text{mm}$, and $A=0$. Each trabecular strut has radius $9/(4\pi)\,\text{mm}$ and the centres are $1.9$\,mm apart. The time point at which the two bone interfaces merge is $t_\text{m} \approx 17$\,days from Eq.~\eqref{circle_analytical}. Figure~\ref{fig_circle_merge} shows the evolution of the interface and velocity at different diffusivities. We perform the simulation first for outward tissue deposition (bone formation) during 34\,days (Fig.~\ref{fig_circle_merge}, left column), then reverse the sign of the velocity from this state (corresponding to bone resorption), and continue the simulation for an additional 49\,days (Fig.~\ref{fig_circle_merge}, right column). Bone adapts to mechanical loading~\cite{Frost1987,Lerebours2016b}, so this scenario represents a first period of mechanical overload, resulting in consolidation of the bone struts by bone formation, followed by a period of return to normal mechanical loads, resulting in resorption of the extraneous bone material by bone resorption. Clearly, this reversal does not lead to the same bone strut geometries as during tissue deposition, i.e., tissue resorption is not simply tissue deposition reversed in time, and bone structure depends on the mechanical loading history~\cite{Lerebours2016b}. The mathematical reason for this time irreversibility is the loss of information that occurs when characteristics collide into shock waves. The weak solutions selected by the level-set method are viscous solutions that satisfy an entropic condition which breaks time-reversal symmetry~\cite{Sethian1999,Leveque2002,AliasPhD}.

The evolution of cell number in these simulations is shown in Figure S1 of the supplementary material. Cell number is reasonably well conserved overall, with some fluctuations ($D=0.0001\,\text{mm}/\text{day}$) or loss ($D=0.01\,\text{mm}/\text{day}$, $D=1\,\text{mm}/\text{day}$) around the time the two interfaces merge, at which extremely acute corners develop.

\paragraph{Fusion and fragmentation of several trabecular bone spicules}
\begin{figure*}[t]  
\centering{\includegraphics[trim={0 0 0 0}, width=0.75\textwidth]{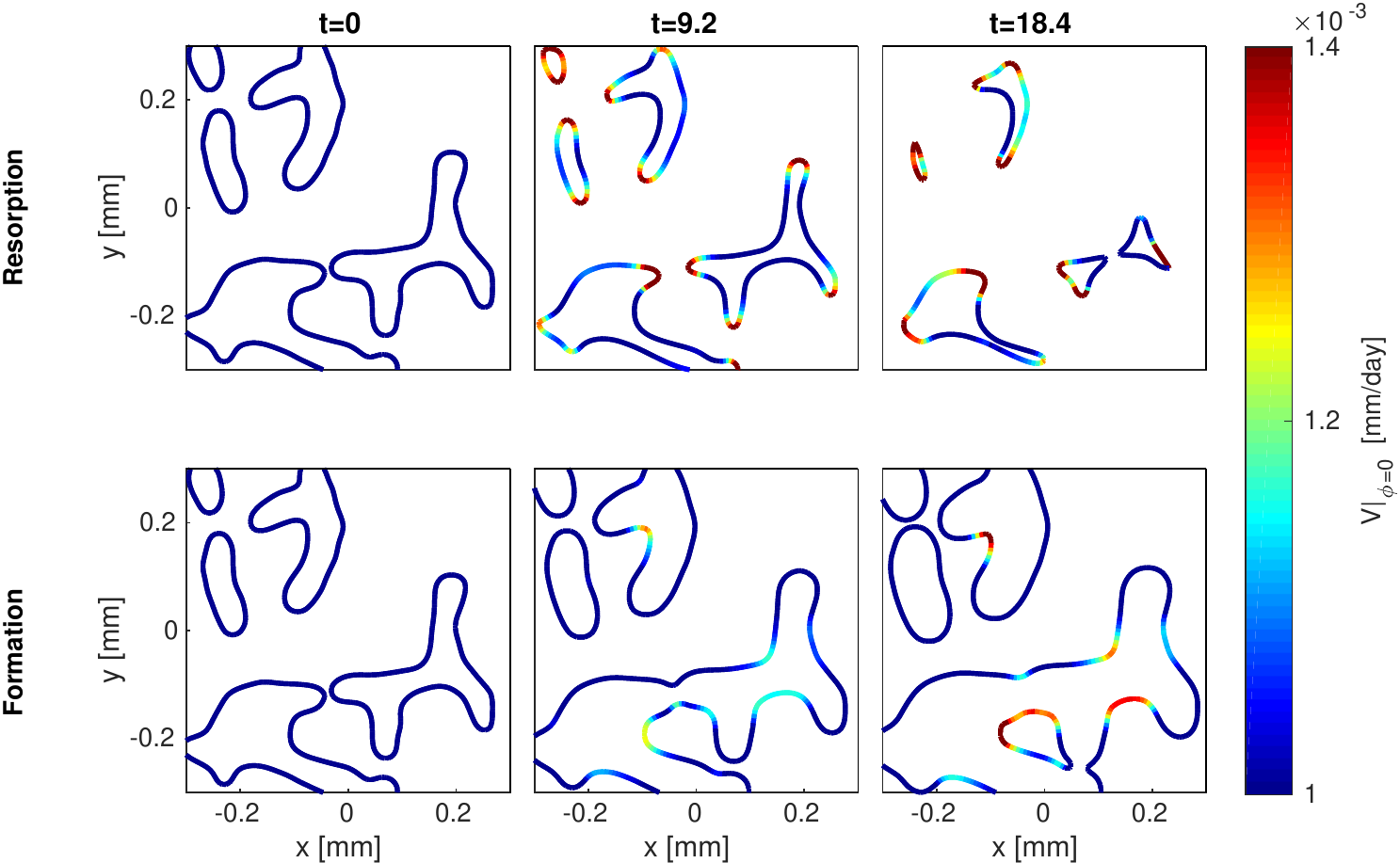}}

\caption{Time snapshots of the concurrent evolution of several trabecular spicules from an experimental image~\cite{Shiga2016} under bone resorption (top) and bone formation (bottom). Simulation parameters: Method~3, $v^0=10^{-3}$mm/day,  $A=0$, $D=0.0001\,\text{mm}^2/\text{day}$, $\epsilon_\text{reinit}=600$, $\Deltaup x=\Deltaup y=0.0086\,\text{mm}$, $\Deltaup t=0.0023\,\text{days}$.}
    \label{fig_complex_trabecular}
\end{figure*}
Figure~\ref{fig_complex_trabecular} shows simulations of bone resorption (top row) and bone apposition (bottom row) from an initial bone interface extracted from a histological section of trabecular bone from Ref.~\cite{Shiga2016}. Bone apposition leads to fusion of initially distinct trabecular spicules, while bone resorption leads to their fragmentation, as would occur for instance during osteoporotic or age-related bone loss~\cite{Maggiano2016b,Bell2001,Kinney1998}. We note that during bone resorption, some trabecular spicules disappear, leading to an inevitable loss of cells that is not entirely due to numerical inaccuracies, see Figure S2 of the supplementary material.

\paragraph{Curvature-controlled tissue growth in bioscaffold}
\begin{figure}[b]
\centering{\includegraphics[width=1\columnwidth]{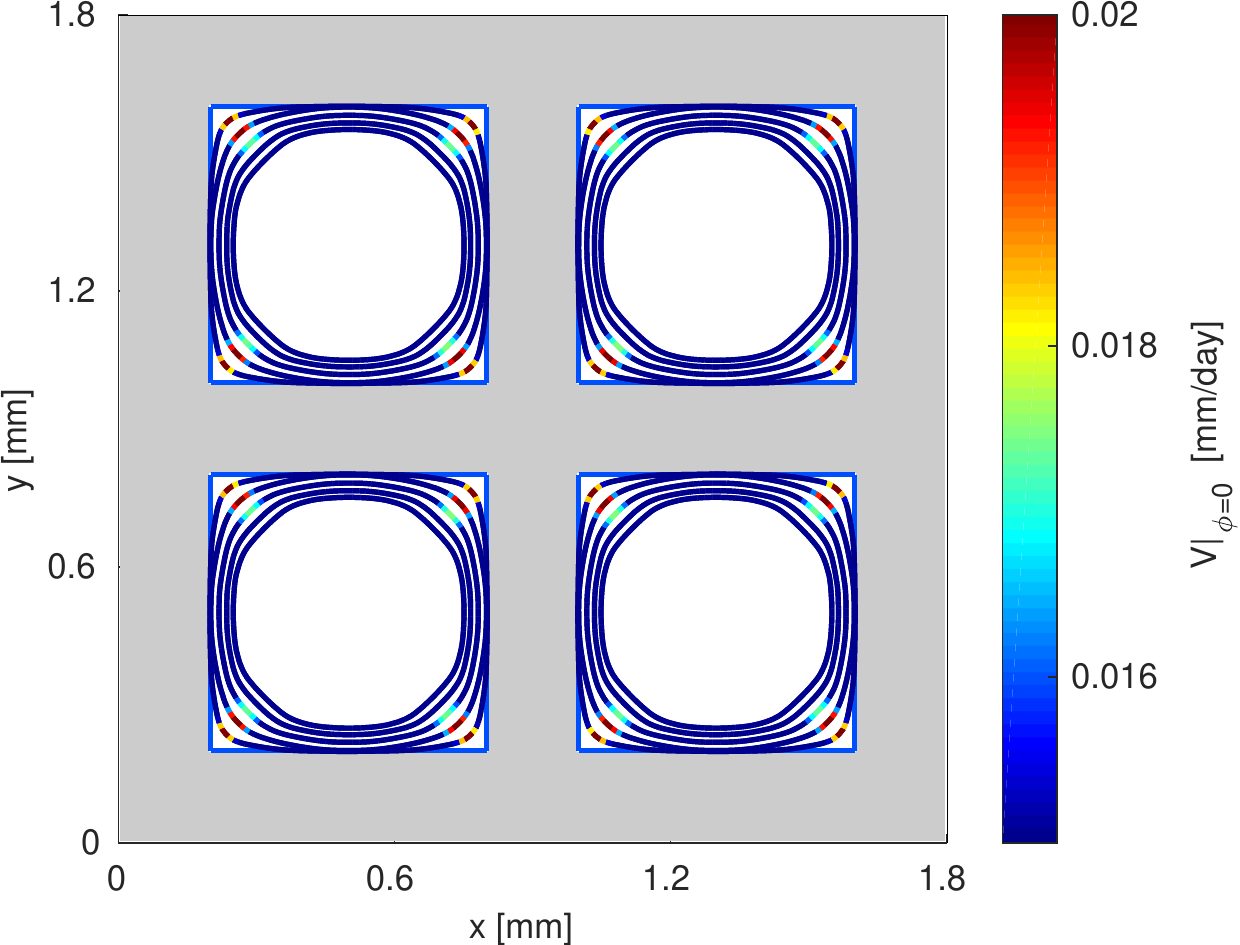}}

\caption{Simultaneous neotissue growth in a bioscaffold consisting of four disconnected pores. The evolution assumes that tissue-synthesing cells are depleting at rate $A=0.1\,\text{day}$ and do not form new tissue on flat or concave areas of the tissue substrate. Simulation parameters: Method 3, $D=0.0001\,\text{mm}^2/\text{day}$, $\epsilon_\text{reinit}=60$, $\Deltaup x = \Deltaup y = 0.017$\,mm and $\Deltaup t = 0.006\,\text{days}$.}
\label{fig_bioscaffold}        
\end{figure}
We now consider an application of the model to the production of neotissue in tissue engineering scaffolds. The model is adapted slightly to prevent tissue formation where the tissue substrate is convex, as suggested experimentally~\cite{Rumpler2008,Bidan2012, Bidan2013, Bidan2013b}. To this effect, we modify the evolution equation of the interface, Eq.~\eqref{phi}, by multiplying the velocity field $V$ with a curvature-dependent Heaviside function:
\begin{align}
\phi_t + H(\kappa) V|\nabla \phi| &= 0,    \label{eqn7_level_set_phi2} 
\end{align} 
where
\begin{align}
H(\kappa) = \left\{ 
\begin{array}{ccc}
1, & \text{ if }  \kappa < 0 &\quad \text{ (concave), }\\
0, & \text{ if }  \kappa \geq 0  & \quad \text{ (flat and convex). }
\end{array}
\right.
\end{align}
Figure~\ref{fig_bioscaffold} represents the simultaneous infilling of four disconnected pores represented by a single level-set function, with a diffusivity $D=0.0001$ mm$^2$/day and a cell depletion rate $A=0.1$/day. The decrease in cell number with time (Supplementary material, Figure S3) matches the theoretical decrease due to cell depletion.

Typically, the porous space of bioscaffolds is extremely complex and would be hard to represent using explicit parameterisations. As illustrated here, our model enables the consideration of cell-specific behaviour, such as curvature-dependent secretory rate (implemented via the function $H(\kappa)$ above), and cell depletion rate $A$. Depletion of active cells is important to explain tissue-deposition slowdown observed experimentally in vitro~\cite{Alias2017} and in-vivo~\cite{Alias2018}. This depletion, and more generally, the inclusion of cell behaviours is not captured by mean curvature flow models of tissue growth~\cite{Bidan2012,Bidan2013,Guyot2014,Guyot2015,Guyot2016,Sanaei2019}.

\paragraph{Three-dimensional tissue growth}

Finally, we consider the growth and shrinkage of a three-dimensional spherical tissue. This geometry allows us to compare numerical solutions with analytic solutions. This example may represent spherical pore infilling or tissue formation atop spherical caps. It may also capture the influence of geometry on the rate of tumour spheroid growth or shrinkage. Avascular tumour spheroids grow by a thin proliferative rim near the surface of the cancerous tissue~\cite{Sutherland1988}. Alternatively, spheroids may shrink, due to, for example, the action of immune cells. In both cases, the normal velocity of the tumour spheroid boundary can be modelled as $v = k \rho$, where the surface cell density $\rho$ is influenced by curvature-induced spreading or crowding (but may also be influenced by other effects). For spheroid growth, the expansion of the tumour is driven by the proliferation of cells in a thin outer layer, which creates new cancerous tissue volume; $\rho$ represents the surface density of the tumour cells in the thin proliferative rim and $k$ corresponds to the tissue volume produced by a proliferative cell per unit time (including daughter cells). As the tumour size increases, more tissue is required to increase the spheroid radius (spreading effect). For tumour shrinkage, $\rho$ represents the surface density of immune cells, such as macrophages, and $k$ corresponds to the volume of tissue they each resorb per unit time. As the tumour shrinks, each immune cell has less surface area of the tumour to resorb (crowding effect on $\rho$). Assuming rotation symmetry, constant proliferation rate~$P$ and constant cell depletion rate~$A$, the evolution of the radius of the spheroid is given by
\begin{align}\label{R-analytic}
    R(t) = R_0 \left[1+\frac{3 v^0}{R^0}\left(\frac{\e^{(P-A)t}}{P-A}-1\right)\right]^{1/3},
\end{align} 
where $v^0=\pm k\rho^0$ is the initial speed of the spheroid boundary (positive for growth, negative for shrinkage) and $R^0$ is the initial radius. This expression can be found by setting $\rho = \frac{N}{4\pi R^2}$, where $N=N^0 \e^{(P-A)t}$ is the total number of cells and $4\pi R^2$ is the surface area of the spheroid, and then solving the single ordinary differential equation $R_t = v = k \rho = \frac{k N^0\e^{(P-A)t}}{4\pi R^2}$ for $R(t)$.

\begin{figure*}[t!]	
	\centering{\includegraphics[trim={0 50 0 0}, width=0.95\textwidth]{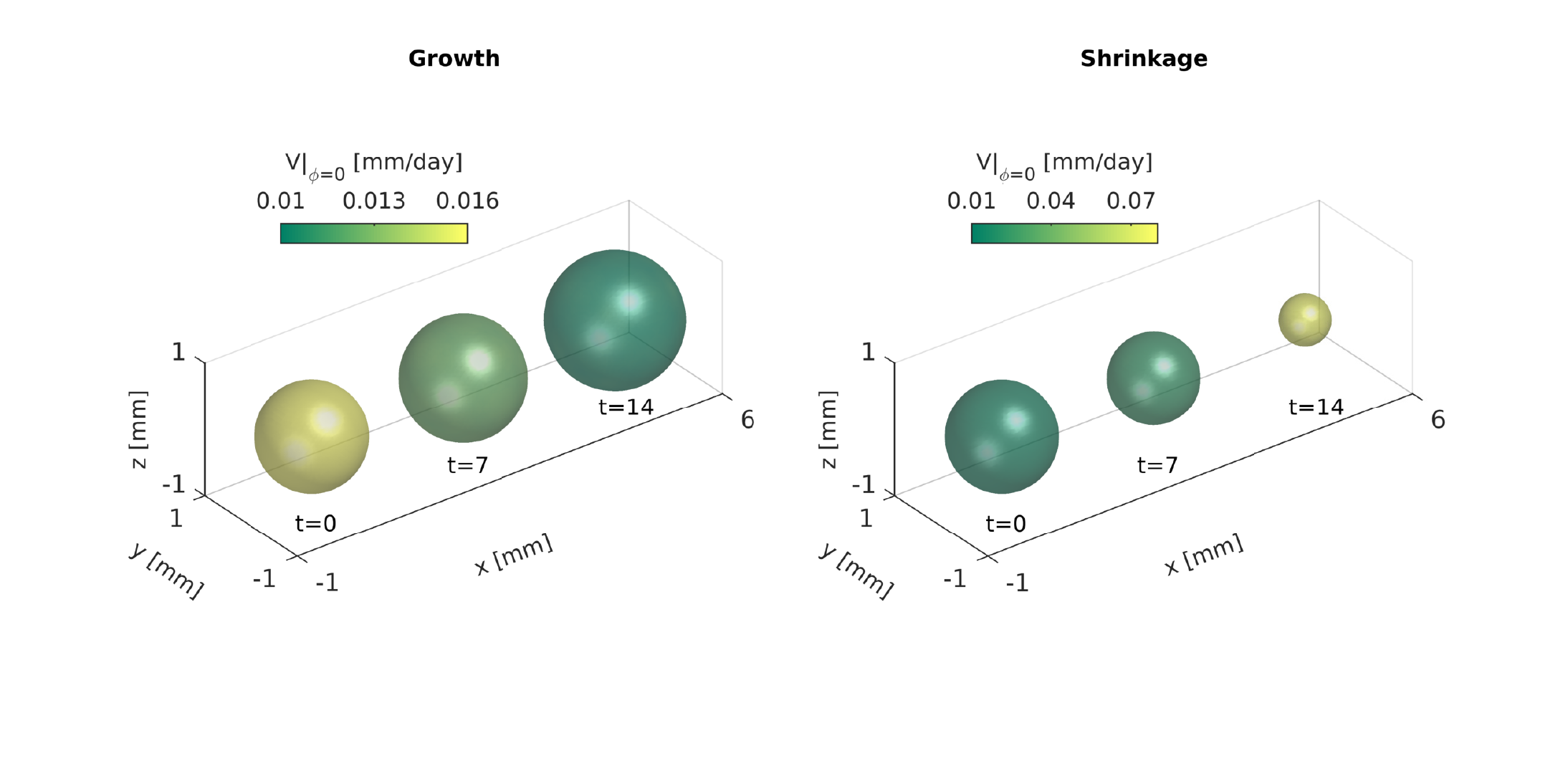}}
	\caption{Tumour spheroid growth (left) and shrinkage (right) in three dimensions in the hyperbolic curvature flow model Simulation parameters: Method 3, $D=0.0001\,\text{mm}^2/\text{day}$, $P-A=0$, $v^0=0.016$\,mm/day, $R^0=0.75$\,mm, $\epsilon_\text{reinit}=300$, $\Deltaup x = \Deltaup y = 0.05$\,mm and $\Deltaup t = 0.028\,\text{days}$.}
	\label{fig_3d}
\end{figure*}

\begin{figure}[b]
	\centering{\includegraphics[width=0.9\columnwidth]{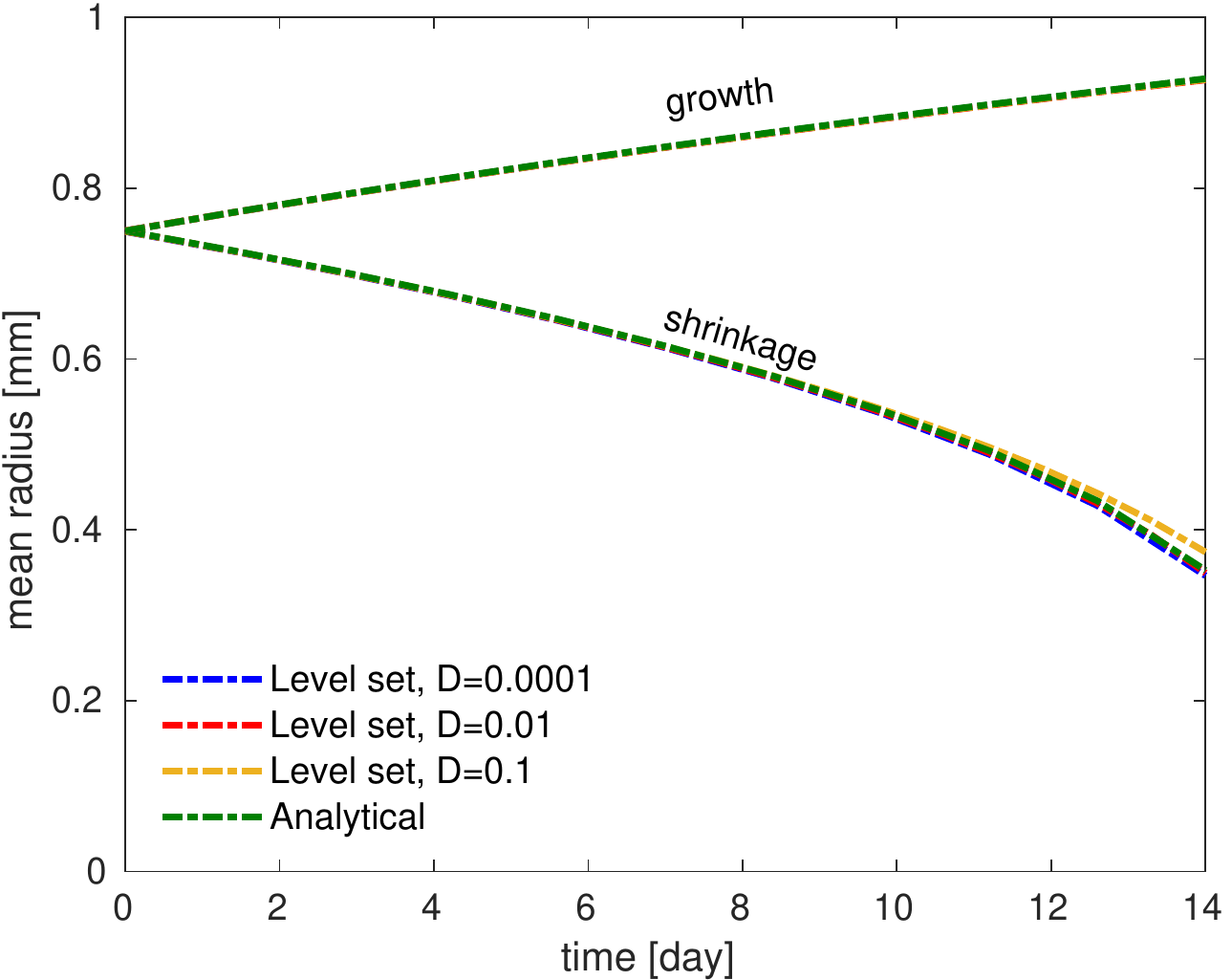}}
	\caption{Comparison between analytic (green) and numeric solutions of growing and shrinking spheroid tissues with $P-A=0$ at $D=0.0001\,\text{mm}^2/\text{day}$ (blue), $D=0.01\,\text{mm}^2/\text{day}$ (red), and $D=0.1\,\text{mm}^2/\text{day}$ (yellow). Same simulation parameters as in Figure~\ref{fig_3d} except $\Delta t=0.014$\,days for $D=0.01\,\text{mm}^2/\text{day}$ and $D=0.1\,\text{mm}^2/\text{day}$.}
	\label{fig_3d_cellnum}
\end{figure}

Figure~\ref{fig_3d} shows time snapshots of the time evolution of a spheroid during growth (left) and shrinkage (right) when proliferation is balanced by cell depletion ($P=A$), so that total cell number is conserved. By symmetry, cell diffusion is irrelevant, but still included in the numerical solutions. There is excellent agreement between the analytic solution given by Eq.~\eqref{R-analytic} and the evolution of the radius of the numerical solutions during growth and shrinkage for a wide range of diffusivities (Figure~\ref{fig_3d_cellnum}).

The growth of tumour spheroids is of course more complex than the simple model above. Proliferation and depletion likely depend on time, and on other factors not considered here, such as nutrient and oxygen intake~\cite{Sutherland1988}. The aim of this model is only to emphasise the geometric influence of curvature onto growth rate in tumour spheroids. Other influences can be added to this model as appropriate to represent experimental data, but the strength of these additional influences must be adjusted taking into account the mechanistic influence of geometry that we have considered. In other continuum models of tumour spheroid growth, the velocity of the spheroid boundary is often defined by a Darcy flow proportional to the gradient of pressure, where pressure build-up is due to cell proliferation~\cite{Macklin2005,Lowengrub2010,Wise2008}. The hyperbolic curvature flow model performs a similar role to pressure in accounting for increased velocity where there is pressure build-up, i.e., where more tissue is produced per unit time. The two main differences are that (i) the hyperbolic curvature flow model does not represent pressure heterogeneities leading to tissue displacement, it represents directly the spatial redistribution of new tissue volume elements; (ii) dissipation/viscosity is built into the Darcy flow model (it corresponds to an overdamped regime and is thus closer to a mean curvature flow model), whereas in the hyperbolic curvature flow model, dissipation is added explicitly as an additional diffusive driving force, whose strength can lead to different interface movement patterns.

\section{Conclusions}

\label{section6_discussion}
The level-set based method proposed in this paper enables us to simulate the co-evolution of tissue interfaces and tissue-synthesising cells in complex geometries, in two and three dimensions
. These simulations account for the influence of mechanistic spatial constraints imposed by tissue shape onto local tissue growth rates.  Mathematically, this corresponds to solving curvature flows of the hyperbolic type, in which interface velocity is determined dynamically by surface-bound processes. This method requires the introduction of an additional Eulerian field to the level set function, which represents the anticipated value of surface cell density at future locations of the interface. The level set function and cell density field are strongly coupled with each other. This coupling is responsible for the rich set of interface behaviour of hyperbolic curvature flows, which includes oscillatory motion, sideways shock propagation, and interface smoothing~\cite{Alias2017}. Comparison with simulations that use explicit parameterisations shows that these different behaviours are well captured by the level-set method proposed.

Importantly, we find that a good indicator of numerical accuracy of the method is provided by tracking the total number of tissue-synthesising cells along the interface with time. Generally, numerical nonconservation is increased at developing cusps in the interface, but reinitialisation of the level-set function and reinitialisation of the velocity field by orthogonal extrapolation help minimise nonconservation in most cases. At very low diffusivities, however, simulations were found to be more accurate without these re-initialisations.

The main advantage of this level-set method for hyperbolic curvature flows is to allow simulations of complex evolving topological situations, that include fragmentation of the interface, and fusion of initially distinct regions of the interface. We have applied the method to the simulation of biological tissue growth to several such complex geometric situations, greatly extending the applicability of such flows to real situations.

Finally, the numerical algorithms we have used may be improved, particularly at developing cusps and where interface regions merge, by using more advanced estimations of unit normals and curvature~\cite{Ervik2014,Coquerelle2016,Vogl2016}. Particle level set methods~\cite{Enright2002} and conservative level set methods~\cite{Olsson2005,Olsson2007} have been developed to help preserve mass in simulations of multi-phase fluid flows. While these techniques help preserve volumetric fluid mass, it is possible that similar techniques could also help preserve interfacial mass. 

Further improvements to the model may include competition for space between two biological tissues. The method could be adapted to this situation by including the mechanical and biochemical interactions between the two tissues into the growth dynamics. Current mathematical models explore such interactions in one dimension~\cite{Murphy2019,Lorenzi2019}.

\subsubsection*{Declarations of interest}

\footnotesize
None
\normalsize

\subsubsection*{Author contributions}

\footnotesize
MAA and PRB conceived and designed the study; MAA performed the numerical simulations; MAA and PRB analysed the data; MAA drafted the article; and both authors edited the article and gave final approval for
publication.
\normalsize

\subsubsection*{Acknowledgments}


\footnotesize
MAA is supported by the Geran Galakan Penyelidik Muda from the Universiti Kebangsaan Malaysia (grant number GGPM-2018-067). PRB gratefully acknowledges the Australian Research Council for Discovery Early Career Research Fellowship (project No.~DE130101191).
\normalsize

\bibliographystyle{WileyNJD-AMA}
\bibliography{reference}


\clearpage \newpage

\title{\bf A level-set method for the evolution of cells and tissue during curvature-controlled growth -- Supplementary material}
\author{Mohd Almie Alias, Pascal R Buenzli}

\date{}
\maketitle

\setcounter{figure}{0}
\setcounter{section}{0}	
\renewcommand{\thefigure}{S\arabic{figure}}

\section{Fusion of two trabecular bone struts and time irreversibility}\label{appx:two_circles}

\begin{figure}[!h]
	
	\centerline{
		\includegraphics[trim={30 0 70 0}, width=1.1\linewidth,clip]{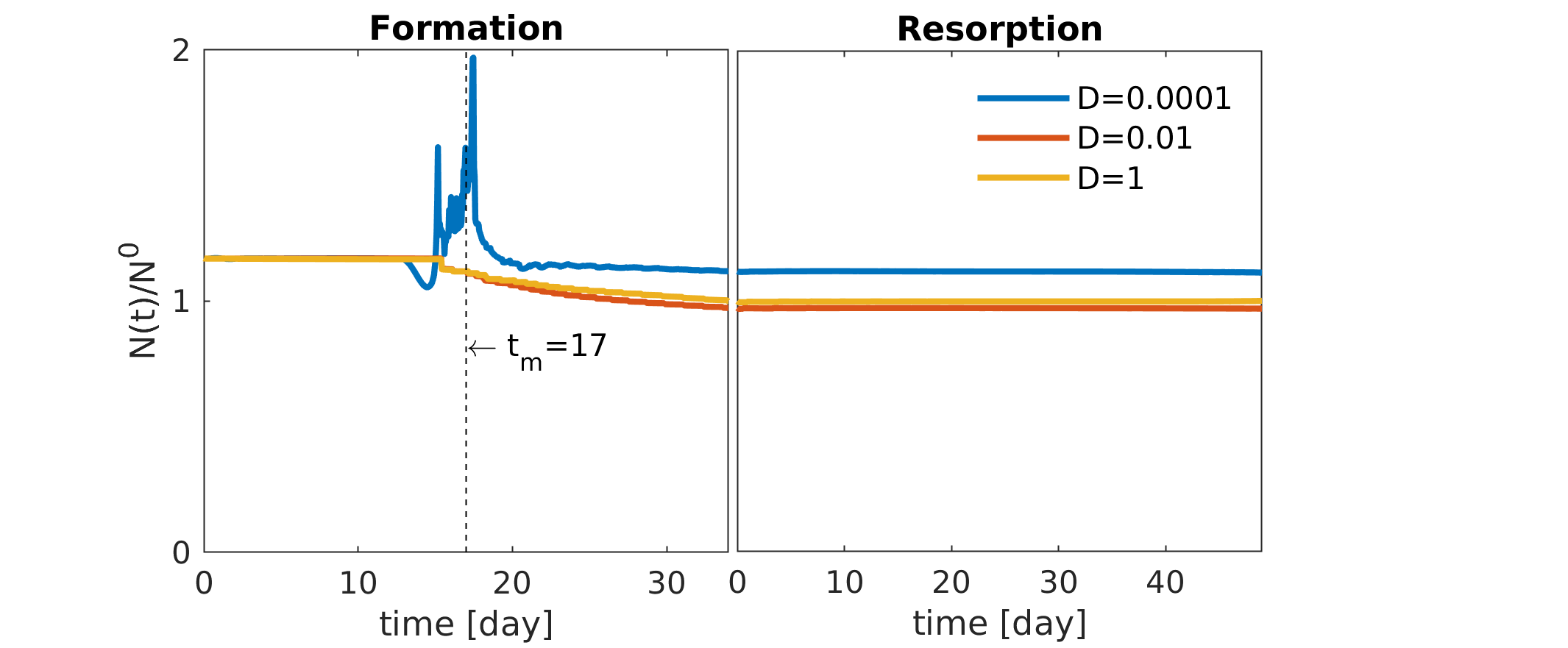}}
	
	\caption{Evolution of normalised cell number for the outward motion during tissue deposition (left) followed by inward motion during tissue resorption of the simulation of two trabecular bone struts in Figure 6. The vertical dashed line corresponds to the time at which the two circles merge.}
	\label{}
\end{figure}

\section{Fusion and fragmentation of trabecular bone spicules}\label{appx:trabecular_spicules}

\begin{figure}[!h]
	
	\centerline{
		\includegraphics[trim={0 0 0 0}, width=0.9\linewidth,clip]{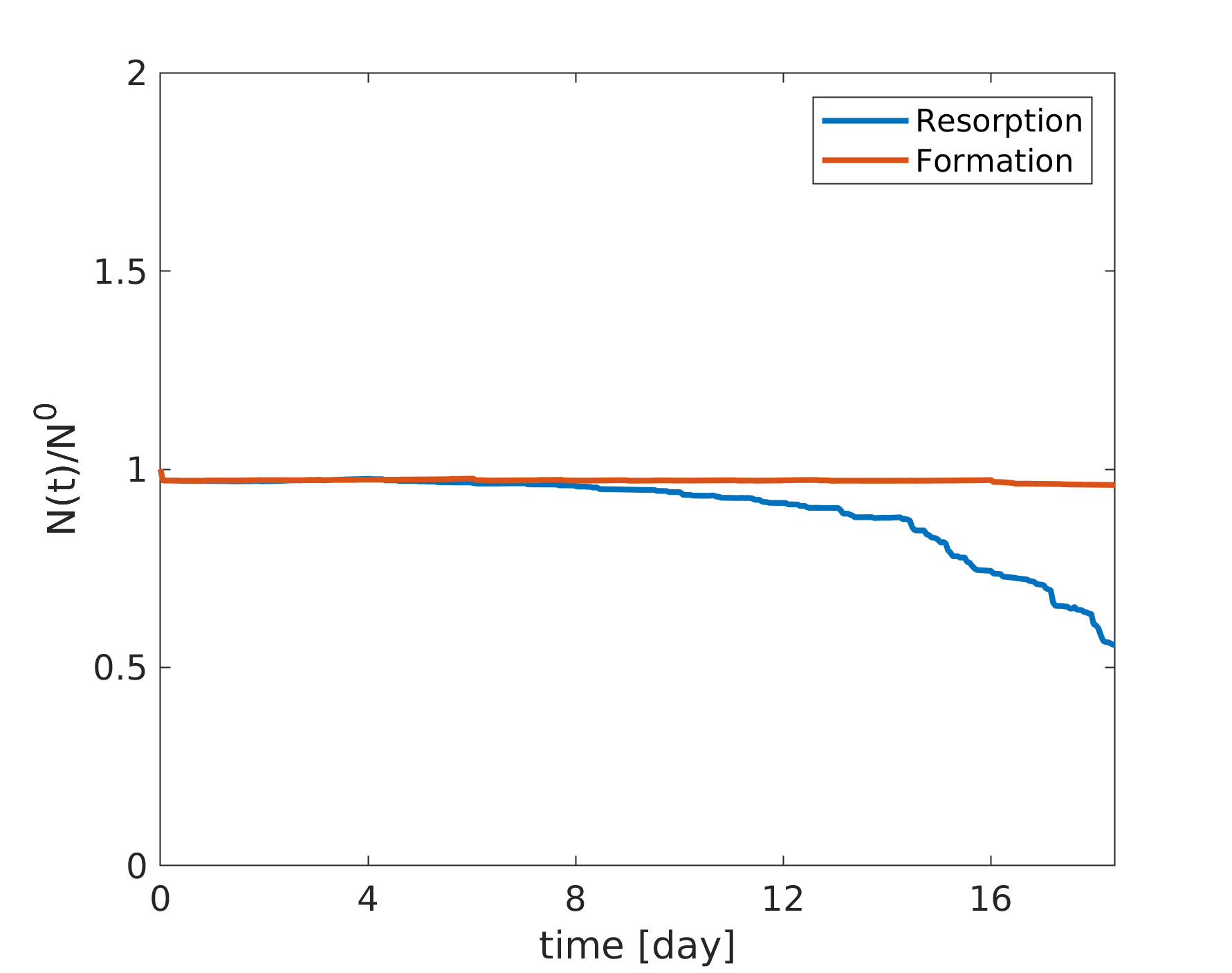}}
	
	\caption{Evolution of normalised cell number for 
		simulations of bone apposition (outward motion -- red) and bone resorption (inward motion -- blue) that causes fusion and fragmentation of the bone in Figure 7.}
	\label{}
\end{figure}

\section{Curvature-controlled tissue growth in bioscaffold}\label{appx:bioscaffold}

\begin{figure}[!h]
	
	\centerline{
		\includegraphics[trim={0 0 0 0}, width=0.9\linewidth,clip]{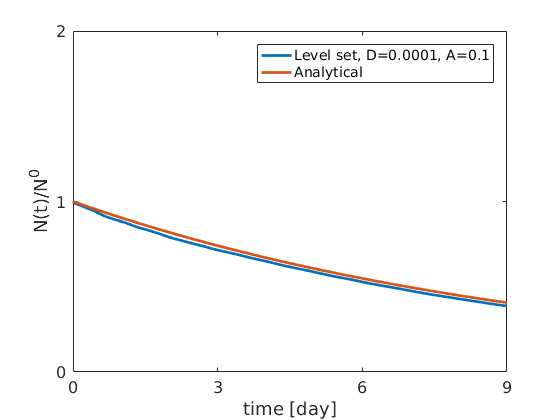}}
	
	\caption{Evolution of normalised cell number during simultaneous infilling of four disconnected pores in a bioscaffold in Figure 8. The analytical normalised cell number is $N(t)/N^0 = \e^{-A t}$.}
	\label{}
\end{figure}

\end{document}